\documentclass[10pt,a4paper]{article}
\usepackage[OT1]{fontenc}
\usepackage[latin1]{inputenc}
\usepackage[english]{babel}
\usepackage{amsfonts}
\usepackage{amsthm}
\usepackage{type1cm} 
\usepackage{graphicx}
\usepackage{setspace}
\usepackage[vcentering,dvips]{geometry} 
\newtheorem{teo}{Theorem}[section]
\newtheorem{prop}[teo]{Proposition}
\newtheorem{lema}[teo]{Lemma}
\newtheorem{coro}[teo]{Corollary}
\newtheorem{defi}[teo]{Definition}
\theoremstyle{definition}
\newtheorem{rem}[teo]{Remark}

\begin{document}

\title{Nonpositive Curvature: a Geometrical Approach to Hilbert-Schmidt Operators}
\date{}
\author{Gabriel Larotonda\footnote{Instituto de Ciencias, Universidad Nacional de General Sarmiento, JM Gutierrez 1150 (1613) Los Polvorines, Buenos Aires, Argentina. e-mail: glaroton@ungs.edu.ar}}

\maketitle

\begin{spacing}{.95}
\abstract{\footnotesize{\noindent We give a Riemannian structure to the set $\Sigma$ of positive invertible unitized Hilbert-Schmidt operators, by means of the trace inner product. This metric makes of $\Sigma$ a nonpositively curved, simply connected and metrically complete Hilbert manifold. The manifold $\Sigma$ is a universal model for symmetric spaces of the noncompact type: any such space can be isometrically embedded into $\Sigma$. We give an intrinsic algebraic characterization of convex closed submanifolds $M$. We study the group of isometries of such submanifolds: we prove that $G_M$, the Banach-Lie group generated by $M$, acts isometrically and transitively on $M$. Moreover, $G_M$ admits a polar decomposition relative to $M$, namely $G_M\simeq M\times K$ as Hilbert manifolds (here $K$ is the isotropy of $p=1$ for the action $I_g: p\mapsto gpg^*$), and also $G_M/K\simeq M$ so $M$ is an homogeneous space. We obtain several decomposition theorems by means of geodesically convex submanifolds $M$. These decompositions are obtained \textit{via} a nonlinear but analytic orthogonal projection $\Pi_M:\Sigma\to M$, a map which is a contraction for the geodesic distance. As a byproduct, we prove the isomorphism $NM\simeq\Sigma$ (here $NM$ stands for the normal bundle of a convex closed submanifold $M$). Writing down the factorizations for fixed ${\rm e}^a$, we obtain ${\rm e}^a={\rm e}^x{\rm e}^v{\rm e}^x$ with ${\rm e}^x\in M$ and $v$ orthogonal to $M$ at $p=1$. As a corollary we obtain decompositions for the full group of invertible elements $G\simeq M\times \exp(T_1M^{\perp})\times K$.}}
\end{spacing}

\section{Introduction}

The aim of this paper is to relate the algebraic and spectral properties of the Banach algebra of unitized Hilbert-Schmidt operators, with the metric and geometrical properties of an underlying manifold $\Sigma$. This is a paper on \textit{applied} nonpositively curved geometry because we first show how the familiar properties of the operator algebra translate into geometrical notions, and then we use the tools of geometry in order to prove new results concerning the operator algebra.

\smallskip

In this paper we study the cone of positive invertible Hilbert-Schmidt operators (extended by the scalar operators) on a separable infinite dimensional Hilbert space $H$. The metric in the tangent space at the identity is given by the trace of the algebra. The local structure induced by the metric is smooth and quadratic; it can be situated in the context of the theory of infinite dimensional Riemann-Hilbert manifolds of nonpositive curvature (\textit{cf.}  Cartan-Hadamard manifolds, as introduced by Lang \cite{lang0}, McAlpin \cite{mcalpin}, Grossman \cite{grossman1} and others). It is then a paper on Riemannian geometry. On the other hand, since the manifold $\Sigma$ is clearly not locally compact, some of the standard results for Hadamard manifolds require a different approach. The geometry is then related to the geometry of the metric spaces in the sense of Aleksandrov \cite{ball0}. It turns out that the notion of convexity (together with the fact that $\Sigma$ is a simply connected and globally nonpositively curved geodesic length space) plays a key role in our constructions. It is then a paper on \textit{metric} geometry.

\smallskip

Through the years, several authors have studied the relationship of geometry and algebra in sets of positive operators, with different approaches that led to a variety of results. In his 1955's paper \cite{mostow}, G.D. Mostow gave a Riemannian structure to the set $M_n^+$ of positive invertible matrices; the induced metric makes  of $M_n^+$ a nonpositively curved symmetric space. Mostow showed that the algebraic concept behind the notion of convexity is that of a \emph{Lie triple system}, which is basically the real part of a given involutive Lie algebra $\mathfrak g$. The geometry of bounded positive operators in an infinite dimensional Hilbert space was studied by G. Corach, H. Porta and L. Recht \cite{cpr6}\cite{cpr2}\cite{pr1} among others, using functional analysis techniques. This area of research is currently very active (see \cite{cm1}\cite{cm2} for a list of references).

\subsection{Main results}

In this paper we study the geometry of a Hilbert manifold $\Sigma$ which is modeled on the operator algebra ${\cal H}_{\mathbb C}$ of unitized Hilbert-Schmidt operators. In Section \ref{intro} we introduce the objects involved and prove some elementary results. The manifold $\Sigma$ is the set of positive invertible operators of ${\cal H}_{\mathbb C}$. Let ${   {\cal H}_{\mathbb C}  }^{\bullet}$ be the classical Banach-Lie group of invertible (unitized) Hilbert-Schmidt operators \cite{pharpe}. The manifold $\Sigma$ has a natural ${   {\cal H}_{\mathbb C}  }^{\bullet}$-invariant metric $<x,y>_{_p}=<xp^{-1},p^{-1}y>_{_2}$, which makes it nonpositively curved (we define $<\alpha+a,\beta+b>_{_2}= \alpha\overline{\beta}+4tr(b^*a)$ whenever $\alpha,\beta\in\mathbb C$ and $a,b$ are Hilbert-Schmidt operators). Let $\mbox{e}^x$ and $\exp(x)$ stand for the usual analytic exponential, \textit{i.e.} $\exp(x)=\sum_{n\ge 0}\frac{x^n}{n!}$. This map is injective when restricted to ${\cal H}_{\mathbb R}$, the set of self-adjoint operators. Let $\ln(p)$ stand for its real analytic inverse. We have $\exp({\cal H}_{\mathbb R})=\Sigma\subset {\cal H}_{\mathbb R}$, and the exponential map induces a diffeomorphism onto its image, so we identify the tangent space at any point of the manifold $\Sigma$ with the set of self-adjoint operators ${\cal H}_{\mathbb R}$, namely $T_p\Sigma\simeq {\cal H}_{\mathbb R}$ for any $p\in\Sigma$. In Section \ref{desiguals} we prove

\medskip

{\bf Theorem A: }
{\it For $p,q\in\Sigma$, the geodesic obtained from Euler's equation by solving Dirichlet's problem is the smooth curve $\gamma_{pq}(t)=p^{\frac12}(p^{-\frac12}qp^{-\frac12} )^tp^{\frac12}$, hence 
$$
{\rm Exp}_p(v)=p^{\frac12}\exp(p^{-\frac12}\;v\;p^{-\frac12} )p^{\frac12}
$$
is the Riemannian exponential of $\Sigma$, for any $v\in T_p\Sigma$. Both ${\rm Exp}_p:T_p\Sigma\to \Sigma$ and its differential map $d\left({\rm Exp}_p\right)_v:T_p\Sigma\to T_{{\rm Exp}_p(v)}\Sigma$ are $C^{\omega}$-isomorphisms for any $p\in\Sigma$, $v\in T_p\Sigma$. The curve $\gamma_{pq}$ is the shortest piecewise smooth path joining $p$ to $q$, hence 
$$
{\rm dist}(p,q)=\|\ln\bigl(p^{\frac12}q^{-1}p^{\frac12}\bigr)\|_{_2}
$$
is the distance in $\Sigma$ induced by the Riemannian metric. The metric space $\left(\Sigma,{\rm dist}\right)$ is complete, and it is globally nonpositively curved.}

\medskip

The curve obtained \textit{via} Calderón's method of complex interpolation \cite{caldera} between the quadratic norms $\|\cdot\|_{_p}$ and $\|\cdot\|_{_q}$ is exactly the short geodesic in $\Sigma$ joining $p$ to $q$ (the proof of \cite{acms} can be adapted almost verbatim).

\medskip

In \cite{grossman1}, N. Grossman proves that the inequality 
\begin{equation}\label{delta}
\|d\left({\rm Exp}_p\right)_v(w)\|_{_p}\ge \|w\|_{_{{\rm Exp}_p(v)}}
\end{equation}
leads to the minimality of geodesics in a simply connected, complete Hilbert manifold. This approach is also carried out by McAlpin \cite{mcalpin}. The following operator inequality involving the differential of the usual exponential map
\begin{equation}\label{dddd}
\|{\rm e}^{-x/2}\;d\exp_x(y) {\rm e}^{-x/2}\|_{_2}\ge \|y\|_{_2}
\end{equation}
is the translation to our context of the inequality (\ref{delta}) above. The convexity of Jacobi fields can be deduced from the non positiveness of the sectional curvature, hence the proof of eqn. (\ref{dddd}) stems in our context from the Cauchy-Schwarz inequality for the trace inner product. We follow the exposition of Lang \cite{lang} on this subject. On the other hand, (\ref{dddd}) can be proved with a direct computation \cite{bhatia}. With this approach the metric completeness of the tangent spaces is not relevant: in Theorem 3.1 of \cite{gandru}, the authors prove the minimizing property of the geodesics in a non complete manifold. The inequality above, in our context, can be also interpreted as the Hyperbolic Cosine Law (see Corollary \ref{triangulines})
$$
a^2\ge b^2+c^2-2bc\cos(\alpha).
$$
Here $a,b,c$ are the lenghts of the sides of any geodesic triangle in $\Sigma$, and $\alpha$ is the angle opposite to $a$. From this inequality also follows that the sum of the inner angles of any geodesic triangle in $\Sigma$ is bounded by $\pi$. 

\medskip

If $A$ is a set of operators, we use $A^+$ to denote the set of positive operators of $A$; note that $( {{\cal H}_{\mathbb C}   }^{\bullet})^+ =\Sigma$. In Section \ref{themain} we show that a submanifold $M\subset\Sigma$ is geodesically convex if and only if its tangent space at the identity $\mathfrak m$ is a Lie triple system. Clearly any such submanifold is nonpositively curved, and Theorem \ref{homo} states:

\medskip

{\bf Theorem B: }
{\it For any geodesically convex, closed submanifold $M=\exp(\mathfrak m)\subset\Sigma$ there exists a connected Banach-Lie group $G_M=\left\langle \exp(\mathfrak m\oplus[\mathfrak m, \mathfrak m])\right\rangle\subset {   {\cal H}_{\mathbb C}  }^{\bullet}$ which acts isometrically and transitively on $M$. Moreover, the polar decomposition of the elements of $G_M$  reduces to $M$ in the sense that $G_M^+=M$. Let $K$ be the isotropy of $1$ for the action; then $K$ is a connected Banach-Lie subgroup of $G_M$ and there is an isomorphism $G_M\simeq M\times K$. In particular any convex submanifold $M$ of $\Sigma$ is an homogeneous space for a suitable Banach-Lie group, which is an analytical subgroup of ${   {\cal H}_{\mathbb C}  }^{\bullet}$. The submanifold $M$ is flat if and only if $M\equiv G_M$ is an abelian Banach-Lie subgroup of ${   {\cal H}_{\mathbb C}  }^{\bullet}$.}

\medskip

The existence of smooth polar decompositions for the involutive Banach-Lie groups can be obtained from the general results of Neeb (\cite{neeb}, Theorem 5.1). Neeb introduces the notion of seminegative curvature (SNC) on Banach-Finsler manifolds $M$, given by the condition of inequality (\ref{delta}) above, plus the condition that $d\left({\rm Exp}_p\right)_v$ should be invertible for any $v\in T_pM$ (the metric of $M$ sould be invariant under parallel transport along geodesics). Neeb proves (Theorem 1.10 of \cite{neeb}) that in a connected, geodesically complete manifold with SNC, the exponential map ${\rm Exp}_p:T_pM\to M$ is a covering map and $M$ is metrically complete, a result which extends that  of Grossman and McAlpin mentioned above to the Banach-Finsler context.

\medskip

The manifold $\Sigma$ can be decomposed by means of any convex closed submanifold $M$. Let $NM$ be the normal bundle of $M$. In Section 5 we prove

\medskip

{\bf Theorem C: }
{\it For any convex closed submanifold $M\subset\Sigma$ there is a nonlinear, real analytic projection $\Pi_M:\Sigma\to M$, which is $\Pi_M$ is contractive for the geodesic distance
$$
{\rm dist}(\Pi_M(p),\Pi_M(q))\le {\rm dist}(p,q)\qquad\mbox{ for any } p,q\in \Sigma.
$$
The point $\Pi_M(p)$ is the (unique) point of $M$ closest to $p$. It can also be viewed as the unique point in $M$ such that there exists a geodesic through $p$ orthogonal to $M$ at $\Pi_M(p)$. The exponential map  $(p,v)\mapsto {\rm Exp}_p(v)$ induces an analytic Riemannian isomorphism $NM\simeq \Sigma$.}

\medskip

Since $\Pi_M(p)$ is the point in $M$ closest to $p$, one can prove the existence of such a point using a metric argument valid in any nonpositively curved geodesic length space \cite{jost}. We choose to give a differential-geometry argument here.

\medskip

In Section 6 we exhibit a decomposition for the submanifold $M=\Delta$ of positive diagonal operators, which is a maximal abelian subalgebra of ${\cal H}_{\mathbb C}$. This decomposition theorem (Theorem \ref{main2}) takes the form of a factorization ${\rm e}^a=d{\rm e}^vd$, where $v$ has null diagonal and $d$ is an invertible diagonal operator. We stress that there is no known algorithm that allows to compute $d$ explicitly (not even if we reduce the problem to $3\times 3$ matrices, that can be thought of as a particular case of the general theory). As a corollary to the decomposition theorems  we obtain

\medskip

{\bf Theorem D: }
{\it Any invertible operator $g\in {   {\cal H}_{\mathbb C}  }^{\bullet}$ admits a unique polar decomposition relative to a fixed closed convex submanifold $M={\rm exp}(\mathfrak m)$. Namely $g={\rm e}^x {\rm e}^v u$ where $x\in \mathfrak m$, $v\in \mathfrak m^{\perp}$ and $u\in {\cal U}({\cal H}_{\mathbb C})$ is a unitary operator. The map $g\mapsto ({\rm e}^x,{\rm e}^v,u)$ is an analytic bijection which gives the isomorphism}
$$
{   {\cal H}_{\mathbb C}  }^{\bullet}\simeq M\times {\rm exp}(\mathfrak m^{\perp})\times {\cal U}({\cal H}_{\mathbb C}).
$$
This isomorphism generalizes the decomposition of $M_n^+$ given in \cite{cpr7}. 

\medskip

In Section 7 we show that the manifold $\Sigma$ can be decomposed by means of a foliation $\{\Sigma_{\lambda}\}_{\lambda>0}$ of totally geodesic submanifolds, namely
$$
\Sigma=\dot{\mathop{\cup}_{\lambda>0}}\Sigma_{\lambda}=\dot{\mathop{\cup}_{\lambda>0}}\{a+\lambda\in \Sigma,\; a=a^* \mbox{ a Hilbert Schmidt operator} \}.
$$
There is a Riemannian isomorphism $\Sigma\simeq\Sigma_1\times \mathbb R_{>0}$ induced by the projection $\Pi_{\Sigma_1}$ of Theorem C above. As an application, we show a decompositon relative to the algebra $M_n^+$ of positive invertible $n\times n$ matrices: fix an $n$-dimensional subspace $S\subset H$, let $P_{_S}$ be the orthogonal projection to $S$ and $Q_{_S}=1-P_{_S}$ the orthogonal projection to $S^{\perp}$. Let $B(S)$ stand for the algebra of bounded linear operators of $S$. Let $R\in {\rm B}(S)^+\simeq M_n^+$, and consider the set 
$$
\mathfrak v=\left\{\left( \begin{array}{lc} 0 & Y^* \\  Y & X \end{array} \right): X=X^*\in {\rm B}(S^{\perp}) \mbox{ a Hilbert-Schmidt operator }, \; Y\in {\rm B}(S,S^{\perp})\right\}.
$$
Let $\cal U({\cal H}_{\mathbb C})$ be the Banach-Lie subgroup of unitary operators in  ${{\cal H}_{\mathbb C} }^{\bullet}$.

\medskip

{\bf Theorem E: }
{\it
For any $g\in {   {\cal H}_{\mathbb C}  }^{\bullet}$ there is a unique factorization $g=\lambda r {\rm e}^v u$ where $\lambda\in\mathbb R_{>0}$, $u\in \cal U({\cal H}_{\mathbb C})$ is a unitary operator, $r=R P_S+ Q_S$ and $v\in \mathfrak v$. In particular
$$
{   {\cal H}_{\mathbb C}  }^{\bullet}\simeq M_n^+\times \exp(\mathfrak v) \times \mathbb R_{>0} \times \cal U({\cal H}_{\mathbb C}).
$$}
The manifold $\Sigma$ can be regarded as a universal model for the symmetric spaces of the noncompact type, namely

\medskip

{\bf Theorem F: }
{\it For any finite dimensional real symmetric manifold $M$ of the noncompact type (\textit{i.e.} with no Euclidean de Rham factor, simply connected and with nonpositive sectional curvature), there is an embedding $M \hookrightarrow \Sigma$ which is a diffeomorphism between $M$ and a closed geodesically convex submanifold of $\Sigma$. If we pull back the inner product on $\Sigma$ to $M$, this inner product is a positive constant multiple of the inner product of $M$ on each irreducible de Rham factor}.

\medskip

The proof of the theorem is straightforward fixing an orthonormal basis of $H$ (see Section \ref{finite}) and recalling the well known result \cite{ebe2} that  for any such space $M$ there is an almost isometric embedding of $M$ into $GL({\mathfrak g})^+$, where ${\mathfrak g}$ is the Lie algebra of the Lie group $I_0(M)$ (the connected component of the identity of the group of isometries of $M$).

\section{Background and definitions}\label{intro}

Let ${\rm B}(H)$ be the set of bounded operators acting on a complex, infinite dimensional and separable Hilbert space $H$, and let ${\sf HS}$ be the bilateral ideal of Hilbert-Schmidt operators of ${\rm B}(H)$. Recall that ${{\sf HS}}$ is a Banach algebra (without unit) when given the norm $\|a\|_{2}=tr(a^*a)^{\frac12}$ (see \cite{simon} for a detailed exposition on trace-class ideals). We will use ${{\sf HS}}^h$ to denote the closed subspace of self-adjoint Hilbert-Schmidt operators. In ${\rm B}(H)$ we define 
$$
{\cal H}_{\mathbb  C}=\{ a+\lambda: \; a\in {{\sf HS}},\; \lambda\in\mathbb C \},
$$ 
the complex linear subalgebra consisting of Hilbert-Schmidt perturbations of scalar multiples of the identity (the closure of this algebra in the operator norm is the set of compact perturbations of scalar multiples of the identity). There is a natural Hilbert space structure for this subspace (where scalar operators are orthogonal to Hilbert-Schmidt operators) which is given by the inner product
$$
<a+\lambda,b+\beta>_{_2}=4tr(ab^*)+\lambda{\overline \beta}.
$$
The algebra ${\cal H}_{\mathbb C}$ is complete with this norm. The model space that we are interested in is the real part of ${\cal H}_{\mathbb C}$, 
$$
{\cal H}_{\mathbb R}=\{ a+\lambda: \;a^*=a,\; a\in {{\sf HS}},\; \lambda\in\mathbb R \},
$$
which inherits the structure of (real) Banach space, and with the same inner product, becomes a real Hilbert space.

\begin{rem}\label{darvuelta}
By virtue of trace properties,  $<xy,y^*x^*>_{_2}=<yx,x^*y^*>_{_2}$ for any $x,y\in {\cal H}_{\mathbb C}$, and  also $<zx,yz>_{_2}=<xz,zy>_{_2}$ for $x,y\in {\cal H}_{\mathbb C}$ and $z\in {\cal H}_{\mathbb R}$.
\end{rem}

Let $\Sigma:=\{A>0:\; A\in {\cal H}_{\mathbb R}\}$ be the subset of positive invertible operators in ${\cal H}_{\mathbb R}$. It is clear that $\Sigma$ is an open set of ${\cal H}_{\mathbb R}$ (for instance, using the lower semi continuity of the spectrum). 

\medskip

\begin{rem}\label{met}
For $p\in\Sigma$, we identify $T_p \Sigma$ with  $ {\cal H}_{\mathbb R}$, and endow this manifold with a (real) Riemannian metric by means of the formula
$$
<x,y>_{_p}:=<p^{-1}x,yp^{-1}>_{_2}=<xp^{-1},p^{-1}y>_{_2}.
$$
Throughout, let  $\|x\|_{_p}:=<x,x>_{_p}^{\frac12}$.  Equivalently, $\|x\|_{_p}=\|p^{-\frac12}xp^{-\frac12}\|_{_2}$.
\end{rem}

\begin{lema}\label{cova}
The covariant derivative in $\Sigma$ (for the metric introduced in Remark \ref{met}) is given by
\begin{equation}\label{covariant}
\left\{\nabla_X Y\right\}_p=\{X(Y)\}_p-\frac12\left( X_p\; p^{-1}\; Y_p+Y_p\;p^{-1}\; X_p  \right).
\end{equation}
Here $X(Y)$ denotes derivation of the vector field $Y$ in the direction of $X$ performed in the linear space ${\cal H}_{\mathbb R}$. 
\end{lema}
\begin{proof}
Note that $\nabla$ is clearly symmetric and verifies all the formal identities of a connection; the proof that it is the Levi-Civita connection relays on the compatibility condition between the connection and the metric, $\frac{d}{dt}<X,Y>_{_{\gamma}}=<\nabla_{\dot\gamma} X,Y>_{_{\gamma}}+<X,\nabla_{\dot\gamma} Y>_{_{\gamma}}$ (see for instance \cite{lang} Chapter VIII, Theorem 4.1). Here $\gamma$ is a smooth curve in $\Sigma$ and $X,Y$ are tangent vector fields along $\gamma$. This identity is straightforward from the definitions and the properties of the trace.
\end{proof}

\smallskip

Let $r^{\alpha}={\rm e}^{\alpha \ln(r)}$ (here $r\in\Sigma,\alpha\in\mathbb R$). The exponential is given by the usual series; note that any positive invertible operator has a real analytic logarithm, which is the inverse of the exponential in the Banach algebra. Note that $aba>0$ whenever $a,b>0$ and also $r^{\alpha}>0$ whenever $r>0$ and $\alpha\in\mathbb R$. 

\medskip

Euler's equation $\nabla_{\dot \gamma}\dot\gamma=0$ for the covariant derivative introduced above reads $\ddot\gamma\ = \dot\gamma\gamma^{^{-1}}\dot\gamma$, 
and it is not hard to see that the (unique) solution of this equation with $\gamma(0)=p$, $\gamma(1)=q$ is given by the smooth curve
\begin{equation}\label{curvas}
\gamma_{pq}(t)=p^{\frac12}( p^{-\frac12}q p^{-\frac12})^t p^{\frac12}.
\end{equation}

\begin{rem}\label{exponencial}
We will use ${\rm Exp}_p:T_p\Sigma\to \Sigma$ to denote the exponential map of $\Sigma$. Differentiating at $t=0$ the curve above, we obtain $\dot{\gamma}_{pq}(0)=p^{\frac12}\ln(p^{-\frac12}q p^{-\frac12})p^{\frac12}$, hence 
$$
{\rm Exp}_p^{-1}(q)=p^{\frac12}\ln(p^{-\frac12}\,q\,p^{-\frac12})p^{\frac12}\quad\mbox{ and  }\quad {\rm Exp}_p(v)=p^{\frac12}\;{\rm exp}(p^{-\frac12}\,v \,p^{-\frac12}) p^{\frac12}.
$$
Note that by the construction above the map ${\rm Exp}_p:T_p\Sigma\to \Sigma$ is surjective (for given $q\in\Sigma$ take $v=p^{\frac12}\ln(p^{-\frac12}\,q\,p^{-\frac12})p^{\frac12}$, then ${\rm Exp}_p(v)=q$). Rearranging the exponential series we get the expressions ${\rm Exp}_p(v)=p\;{\rm e}^{p^{-1}v}={\rm e}^{\,vp^{-1}}p$.
\end{rem}

\begin{lema}\label{invar}
The metric in $\Sigma$ is invariant under the action of the group of invertible elements: if $g$ is an invertible operator in ${\cal H}_{\mathbb C}$, then $I_g(p)=gpg^*$ is an isometry of $\Sigma$.\index{isometry}
\end{lema}
\begin{proof}
First note that for any $\psi\in H$ we have $<gpg^*\psi,\psi>=<pg^*\psi,g^*\psi>=<p\eta,\eta>\;>0$ assuming $p>0$ and $g$ invertible, so $I_g$ maps $\Sigma$ into itself. Also note that $d(I_g)_r(x)=gxg^*$ for any $x\in T_r\Sigma$, hence
$$
\|gxg^*\|_{_{grg^*}}^2=< gxg^*(g^*)^{-1}r^{-1}g^{-1}, (g^*)^{-1}r^{-1}g^{-1}gxg^* >_{_2}=\qquad\qquad\qquad\qquad
$$
$$
\quad = < gxr^{-1}g^{-1}, (g^*)^{-1}r^{-1}xg^* >_{_2}=< xr^{-1}, r^{-1}x >_{_2}=\|x\|_{_r}^2
$$
where the third equality in the above equation follows from Remark \ref{darvuelta}.\end
{proof}

\section{Local and global structure}\label{desiguals}

\subsection{Curvature}\label{ctensor}

We start showing that curvature in this manifold is a measure of noncommutativity, and then give a few definitions, which are necessary because of the infinite dimensional setting. Let $[ \: ,\:]$ stand for the usual commutator of operators, $[x,y]=xy-yx$.

\begin{prop}The curvature tensor for the manifold $\Sigma$ is given by:
\begin{equation}\label{tensor}
{\mathfrak R}_p(x,y)z=-\frac14 \; p\left[  \left[p^{-1}x,p^{-1}y \right]  ,  p^{-1}z   \right].
\end{equation}
\end{prop}
\begin{proof}
This follows from the usual definition ${\mathfrak R}(x,y)=\nabla_x\nabla_y-\nabla_y\nabla_x-\nabla_{[x,y]}$. The formula for $\nabla$ given in Lemma \ref{cova}.
\end{proof}

\begin{defi}
A Riemannian submanifold $M\subset\Sigma$ is \emph{flat} at $p\in M$ if the sectional curvature vanishes for any 2-subspace of $T_pM$. The manifold $M$ is \emph{flat} if it is flat at any $p\in M$. The manifold $M$ is \emph{geodesic} at $p\in M$ if geodesics of the ambient space starting at $p$ with initial velocity in $T_pM$ are also geodesics of $M$. The manifold $M$ is a \emph{totally geodesic} manifold if it is geodesic at any $p\in M$. Equivalently,  $M$ is totally geodesic if any geodesic of $M$ is also a geodesic of  $\Sigma$.
\end{defi}

\begin{prop}
The manifold $\Sigma$ has nonpositive sectional curvature.
\end{prop}
\begin{proof}Let $x,y\in T_p\Sigma$. Let $\overline{x}=p^{-\frac12}xp^{-\frac12}$, $\overline{y}=p^{-\frac12}yp^{-\frac12}$. We may assume that $x,y$ are orthonormal at $p$. A straightforward computation shows that
$$
{\cal S}_p(x,y)=<{\mathfrak R}_p(x,y)y,x>_{_p}=-\frac{1}{4} \left\{<\overline{x}\overline{y}^2,\overline{x}>_{_2}-2<\overline{y}\overline{x}\overline{y},\overline{x}>_{_2}+<\overline{y}^2\overline{x},\overline{x}>_{_2}\right\}.
$$
Since $\overline{x},\overline{y}\in {\cal H}_{\mathbb R}$, $\overline{x}=\lambda+a$ and $\overline{y}=\beta+b$ for $\lambda,\beta\in\mathbb R$ and $a,b\in\sf HS$. The equation reduces to
\begin{equation}\label{curvate}
{\cal S}_p(x,y)=-\frac12\left\{tr(a^2b^2)-tr((ab)^2)\right\}.
\end{equation}
Note that $[z,w]=tr(w^*z)$ is an inner product on ${\sf HS}$, so we have  the Cauchy-Schwarz inequality $tr(w^*z)\le tr^{\frac12}(w^*w)\;tr^{\frac12}(z^*z)
$. Putting $w=ba, z=ab$, we obtain 
$$
tr((ab)^2)=tr(abab)=tr((ba)^*ab)\le tr^{\frac12}(abba)\;tr^{\frac12}(baab)=tr(a^2b^2).\qedhere
$$
\end{proof}

\begin{prop}\label{abe} Let $M\subset \Sigma$ be a submanifold. Assume that $M$ is flat and geodesic at $p\in M$. If $x,y\in T_pM$, then $p^{-\frac12}xp^{-\frac12}$ commutes with $p^{-\frac12}y p^{-\frac12}$.
\end{prop}
\begin{proof}
Since $M$ is geodesic at $p$, the curvature tensor is the restriction of the curvature tensor of $\Sigma$, so in equation (\ref{curvate}) above the right hand term must be zero if $M$ is flat at $p$. But the Cauchy-Schwarz inequality is an equality only if the vectors are linearly dependent; in the notation of the previous theorem, we have $ab=z=\alpha w=\alpha ba$ for some $\alpha\in\mathbb R$; replacing this in the above equation we obtain $\alpha=1$, namely $ab=ba$. Recalling the definitions for $a$ and $b$ we obtain the assertion.
\end{proof}

\subsection{Convexity of Jacobi fields}

Let $J(t)$ be a Jacobi field along a geodesic $\gamma$ of $\Sigma$, \textit{i.e.} $J$ is a solution of the differential equation\index{Jacobi vector field}
\begin{equation}\label{ecuacionjacobi}
D_t^2 J+{\mathfrak R}_\gamma(J,\dot{\gamma})\dot{\gamma}=0
\end{equation}
where $D_t=\nabla_{\dot{\gamma}}$ is the covariant derivative along $\gamma$. We may assume that $J(t)$ is non vanishing, hence
$$
\begin{array}{rl}
\|J\|_{_{\gamma}}^{3}\;\frac{d^2}{dt^2}<J,J>_{_{\gamma}}^{\frac12} = &
- <D_tJ,J>_{_{\gamma}}^2+<J,J>_{_{\gamma}}\;<D_tJ,D_tJ>_{_{\gamma}}\\
& -<J,J>_{_{\gamma}}\;<{\mathfrak R}_{\gamma}(J,\dot{\gamma})\dot{\gamma},J> _{\gamma}.  
\end{array}
$$
The third term is clearly positive and the first two terms add up to a nonnegative number by the Cauchy-Schwarz inequality: $<D_tJ,J>_{_{\gamma}}^2\;\le\; <D_tJ,D_tJ>_{_{\gamma}}\;<J,J>_{_{\gamma}}$. In other words, the smooth function $t\mapsto <J,J>_{_{\gamma}}^{\frac12}=\|J\|_{_{\gamma}}$ is convex, exactly as in the finite dimensional setting.

\subsection{The exponential map}

We present two theorems that, in this infinite dimensional setting, stem from McAlpin's PhD. Thesis (for a proof see \cite{mcalpin} or Theorem 3.7 of Chapter IX in \cite{lang}). First, if one identifies the Riemannian exponential with a suitable Jacobi lift, one obtains
\begin{teo}\label{espandex} The map ${\rm Exp}_p:T_p\Sigma\to\Sigma$ has an expansive differential: 
$$
\|d\left( {\rm Exp}_p\right)_v(w)\|_{_{{\rm Exp}_p(v)}}\ge \| w\|_{_p}.
$$
\end{teo}

\smallskip

This result implies that the differential of the exponential map is injective and has closed range. Playing with the Hilbert structure of the tangent bundle and using the well known identity for operators $Ker(A)^{\perp}=\overline{Ran(A^*)}$, it can be proved that this map is surjective, moreover 

\begin{coro}\label{esiso}
The differential of the Riemannian exponential $d({ \rm Exp}_p)_v:T_p\Sigma\to T_{{\rm Exp}_p(v)}\Sigma$ is a linear isomorphism for any $v\in T_p\Sigma$. Hence, ${\rm Exp_p}:T_pM\to \Sigma$ is a $C^{\omega}$-diffeomorphism.
\end{coro}
The last assertion is due to the fact that the map ${\rm Exp}_p:T_p\Sigma\to \Sigma$ is a bijection (see Remark \ref{exponencial} above).

\subsection{The shortest path and the geodesic distance}\label{geods}

The following inequality is the key to the proof of the fact that geodesics are minimizing. It was proved by R. Bhatia \cite{bhatia} for matrices, and his proof can be translated almost verbatim to the context of operator algebras with a trace, see \cite{gandru}. However since the Riemannian metric in $\Sigma$ is complete, the inequality can be easily deduced from the fact that the norm of a Jacobi field is a convex map (in Theorem \ref{espandex} put $p=1$, $v=x$ and $w=y$):

\begin{coro}\label{de}
If $d\;{\rm exp}_{x}$ denotes the differential at $x$ of the usual exponential map, then for any $x,y\in {\cal H}_{\mathbb R}$
$$
\|d\;{\rm exp}_{x}(y)\|_{_{{\rm e}^{x}}}=\|{\rm e}^{-\frac{x}{2}}d\;{\rm exp}_{x}(y){\rm e}^{-\frac{x}{2}} \|_{_2}\ge \|y\|_{_2}.
$$
\end{coro}

As usual, one measures length of curves in $\Sigma$ using the norms in each tangent space,
\begin{equation}\label{metrica}
L(\alpha)=\int_0^1\|\dot{\alpha}(t)\|_{_{{\alpha}(t)}}\, dt.
\end{equation}
We define the distance between two points $p,q\in \Sigma$ as the infimum of the lengths of piecewise smooth curves in $\Sigma$ joining $p$ to $q$,
$$
{\rm dist}(p,q)=\inf\left\{L(\alpha):\alpha\subset\Sigma, \;\alpha(0)=p, \;\alpha(1)=q\right\}.
$$
Recall (Remark \ref{exponencial} and the paragraph above it) that for any pair of elements $p,q\in \Sigma$, we have the smooth curve $\gamma_{pq}\subset \Sigma,\; \gamma_{pq}(t)=p^{\frac12}(p^{-\frac12}qp^{-\frac12})^t p^{\frac12}
$ joining $p$ to $q$, which is the unique solution of Euler's equation in $\Sigma$. Computing the derivative, we get
$$
\|\dot{\gamma}_{pq}(t)\|_{_{\gamma_{pq}(t)}}\equiv \|\ln(p^{\frac12}q^{-1}p^{\frac12})\|_{_2}=L(\gamma_{pq}).
$$
The minimality of these (unique) geodesics joining two points can be deduced from general considerations \cite{grossman1}, we present here a direct proof.
\begin{teo}\label{postal}
Let $p,q\in \Sigma$. Then the geodesic $\gamma_{pq}$ is the shortest curve joining $p$ and $q$ in $\Sigma$, if the length of curves is measured with the metric defined above (\ref{metrica}).
\end{teo}
\begin{proof}
Let $\alpha$ be a smooth curve in $\Sigma$ with $\alpha(0)=p$ and $\alpha(1)=q$. We must compare the length of $\alpha$ with the length of $\gamma_{pq}$. Since the invertible group acts isometrically, it preserves the lengths of curves. Thus we may act with $p^{-\frac12}$, and suppose that both curves start at $1$, or equivalently that $p=1$. Therefore $\gamma_{1q}(t):=\gamma(t)={\rm e}^{tx}$, with $x=\ln q$. The length of $\gamma$ is then $\|x\|_{_2}$. 
The proof follows easily from the inequality of Corollary \ref{de}. Indeed, since $\alpha$ is a smooth curve in $\Sigma$, it can be written as $\alpha(t)={\rm e}^{\beta(t)}$, with $\beta=\ln \alpha$. Then $\beta$ is a smooth curve of self-adjoint operators with $\beta(0)=0$ and $\beta(1)=x$.
Moreover,
$$
L(\gamma)=\|x\|_{_2}=\|x-0\|_{_2}=\|\int_0^1 \dot{\beta}(t)\;dt\|_{_2}\le \int_0^1 \|\dot{\beta}(t)\|_{_2}\;dt.
$$
On the other hand, by the mentioned inequality, 
$$
\|\dot{\beta}(t)\|_{_2}\le \| {\rm e}^{-\frac{\beta(t)}{2} } \; d\, {\rm exp}_{\beta(t)}(\dot{\beta}(t) ) {\rm e}^{-\frac{\beta(t)}{2}} \|_{_2}=\|  d\, {\rm exp}_{\beta(t)}(\dot{\beta}(t) )\|_{_{{\rm e}^{\beta(t)}}}=\| \dot{\alpha}(t) \|_{_{\alpha(t)}}.\qedhere
$$
\end{proof}

\begin{rem}
The geodesic distance induced by the metric is given by
$$
{\rm dist}(p,q)=\|\ln\bigl(p^{\frac12}\; q^{-1}p^{\frac12}\bigr)\|_{_2}.
$$
\end{rem}
Hence the unique geodesic joining $p$ to $q$ is also the shortest path joining $p$ to $q$. This means that $\left(\Sigma,{\rm dist}\right)$ is a (not locally compact) \textit{geodesic length space} in the sense of Aleksandrov and Gromov \cite{ball0}. These curves look formally equal to the geodesics between positive definite $n\times n$ matrices, when this space is regarded as a symmetric space. 

\begin{coro}\label{con}
If $\gamma$, $\delta$ are geodesics, the map $f:\mathbb R\to \mathbb R$, $t\mapsto {\rm dist}(\gamma(t),\delta(t))$ is convex.
\end{coro}
\begin{proof}
The distance between the points $\gamma(t)$ and $\delta(t)$ is given by the geodesic $\alpha_t(s)$, which is obtained as the $s$ variable ranges in a geodesic square $h(s,t)$ with vertices $\left\{\gamma(t_0),\delta(t_0),\gamma(t_1),\delta(t_1)\right\}$ (the starting and ending points of $\gamma$ and $\delta$). Taking the partial derivative along the direction of $s$ gives a Jacobi field $J(s,t)$ along the geodesic $\beta_s(t)=h(s,t)$ and it also gives the speed of $\alpha_t$. Hence 
$$
f(t)=\int_0^1\|\frac{\partial\alpha_t}{\partial s}(s)\|_{_{\alpha_t(s)}}ds=\int_0^1 \| J(s,t)\|_{_{h(s,t)}}ds.
$$
This equation states that $f(t)$ can be written as the limit of a convex combination of convex functions $u_i(t)= \| J(s_i,t)\|_{_{h(s_i,t)}}$, so $f$ must be convex itself.
\end{proof}

In a recent paper (Corollary 8.7 of \cite{lawson}), the authors prove this property of convexity of the geodesic distance in a general setting concerning nonpositively curved symmetric spaces given by a quotient of Banach-Lie groups.

\begin{lema}\label{este}
For any $x,y\in {\cal H}_{\mathbb R}$ we have
\begin{equation}\label{des}
{\rm dist}({\rm e}^x,{\rm e}^y)=\|\ln({\rm e}^{x/2}{\rm e}^{-y}{\rm e}^{x/2})\|_{_2}\ge \|x-y\|_{_2}
\end{equation}
\end{lema}
\begin{proof}
Take $\gamma(t)={\rm e}^{tx}$, $\delta(t)={\rm e}^{ty}$ and $f$ as in the previous corollary; we may assume that $x,y\in {\sf HS}^h$. Note that $f(0)=0$, hence $f(t)/t\le f(1)$ for any $0<t\le 1$; hence $\lim\limits_{t \to 0^+}f(t)/t\le f(1)$. Now
$$
f(t)/t=\frac{1}{t}\|\ln({\rm e}^{tx/2}{\rm e}^{-ty}{\rm e}^{tx/2})\|_{_2}=tr([\frac{1}{t}\ln({\rm e}^{tx/2}{\rm e}^{-ty}{\rm e}^{tx/2}) ]^2)^{\frac12},
$$
and
$$
\lim\limits_{t\to 0^+}\frac{1}{t}\ln({\rm e}^{tx/2}{\rm e}^{-ty}{\rm e}^{tx/2})=\frac{d}{dt}\mid_{t=0}\ln({\rm e}^{tx/2}{\rm e}^{-ty}{\rm 	e}^{tx/2})=d\ln_1(x-y)=x-y.\qedhere
$$
\end{proof}

\medskip

\begin{coro}\label{triangulines} The inner angles of any geodesic triangle in $\Sigma$ add up to at most $\pi$.
\end{coro}
\begin{proof} Using the invariance of the metric for the action of the group of invertible operators, and squaring both sides of inequality (\ref{des}) in Lemma \ref{este}, we obtain the Hyperbolic Cosine Law: 
\begin{equation}\label{hypcos}
l_i^2\ge l_{i+1}^2+l_{i-1}^2-2l_{i+1}l_{i-1}\cos(\alpha_i).
\end{equation}
Here $l_i$ (i=1,2,3) are the sides of any geodesic triangle and $\alpha_i$ is the angle opposite to $l_i$. These inequalities put together show that one can construct a comparison Euclidean triangle in the affine plane with sides $l_i$. For this triangle with angles $\beta_i$ (opposite to the side $l_i$) we have $l_i^2= l_{i+1}^2+l_{i-1}^2-2l_{i+1}l_{i-1}\cos(\beta_i)
$. This equation together with inequality (\ref{hypcos}) imply that the angle $\beta_i$ is bigger than $\alpha_i$ for $i=1,2,3$. Adding the three angles we have
$\alpha_1+\alpha_2+\alpha_3\le \beta_1+\beta_2+\beta_3=\pi$.
\end{proof}

\begin{prop}\label{compmetri}
The metric space $\left(\Sigma,d\right)$ is complete with the distance induced by the minimizing geodesics.
\end{prop}
\begin{proof}
Consider a Cauchy sequence $\{p_n\}\subset\Sigma$. Again by virtue of inequality (\ref{des}) of Lemma \ref{este}, $x_n=\ln(p_n)$ is a Cauchy sequence in ${\cal H}_{\mathbb R}$. Since Hilbert-Schmidt operators are complete with the trace norm, there is a vector $x\in {\cal H}_{\mathbb R}$ such that $x_n\to x$ in the trace norm. Since the inverse map, the exponential map, the product and the logarithm are all analytic maps with respect  to the trace norm, ${\rm dist}(p_n,{\rm e}^x)=\|\ln({\rm e}^{x/2}  {\rm e}^{-x_n}{\rm e}^{x/2})\|_{_2}\to 0$ when $n\to\infty$.
\end{proof}

\section{Geodesically convex submanifolds}\label{themain}

\begin{defi}A set $M\subset \Sigma$ is \emph{geodesically convex} if for any two given points $p,q\in {M}$, the unique geodesic of $\Sigma$ joining $p$ to $q$ lies entirely in ${M}$. A Riemannian submanifold $M\subset\Sigma$ is complete at $p\in M$ if $\;{\rm Exp}_p^M$ is defined in the whole tangent space and maps onto $M$. The manifold $M$ is  \emph{complete} if it is complete at any point. 
\end{defi}

\begin{rem}The manifold $\Sigma$ is complete; moreover, ${\rm Exp}_p$ is a $C^{\omega}$ (analytic) isomorphism of $\cal H_{\mathbb R}$ with $\Sigma$ for each $p\in\Sigma$. Other notions of completeness are touchy because, as C. J. Atkin shows in \cite{atkin1} and \cite{atkin2}, the Hopf-Rinow Theorem does not necessarily hold in infinite dimensional Banach manifolds.
\end{rem}

These previous notions are strongly related; it is not hard to see that for any Riemannian submanifold $M$ of $\Sigma$, $M$ is geodesically convex if and only if $M$ is complete and totally geodesic. On the other hand, it should be clear from the definitions that whenever $M$ is a convex submanifold of $\Sigma$, $M$ is nonpositively curved.

\subsection{An intrinsic characterization of convexity}\label{convecca}

From now on the term \textit{convex} stands for the longer \textit{geodesically convex}. As before $[\; ,\;]$  denotes the usual commutator of operators in ${\rm B}(H)$. To deal with convex sets the following definition will be useful; assume $\mathfrak m\subset {\cal H}_{\mathbb R}$ is a real linear space.

\begin{defi}
We say that $\mathfrak m$ is a \emph{Lie triple system} if $[[a,b],c]\in\mathfrak m$ for any $a,b,c\in \mathfrak m$. Equivalently, $[x,[x,y]]\in\mathfrak m$ whenever $x,y\in \mathfrak m$.
\end{defi} 

Note that whenever $a,b,c$ are self-adjoint operators, $d=[a,[b,c]]$ is also a self-adjoint operator. So, for any involutive Lie subalgebra of operators $\mathfrak a\subset {\cal H}_{\mathbb C}$ (in particular: for any associative Banach subalgebra), $\mathfrak m=\mathfrak{Re(a)}$ is a Lie triple system in $ {\cal H}_{\mathbb R}$.

\medskip

Assume $M\subset\Sigma$ is a submanifold such that $1\in M$, and $M$ is geodesic at $p=1$. Then $T_1M$ is a Lie triple system, because the curvature tensor at $p=1$ is the restriction to $T_1M$ of the curvature tensor of $\Sigma$, and ${\mathfrak R}_1(x,y)z=-\frac14 [[x,y],z] $. In particular, if $M$ is geodesically convex, $T_1M$ must be a Lie triple system. This weak condition on the tangent space turns out to be strong enough to obtain convexity:

\begin{teo}(Mostow-de la Harpe \cite{mostow}\cite{pharpe})\label{pharpe}
Assume $\mathfrak m\subset {\cal H}_{\mathbb R}$ is a closed subspace, put $M={\rm exp}(\mathfrak m)\subset\Sigma$ with the induced topology and Riemannian metric. Assume further that $\mathfrak m$ is a Lie triple system. Then for any $p,q\in M$ it holds true that $qpq\in M$.
\end{teo}\label{mostro}
\begin{proof}As P. de la Harpe pointed out, the proof of G. D. Mostow for matrices in \cite{mostow} can be translated to Hilbert-Schmidt operators without any modification: we give a sketch of the proof here. Assume $p={\rm e}^x$, $q={\rm e}^y$, and consider the curve ${\rm e}^{\alpha(t)}={\rm e}^{ty}{\rm e}^x{\rm e}^{ty}$. Then it can be proved that $\dot\alpha(t)=G(\alpha(t))$ with $G$ a Lipschitz map that sends $\mathfrak m$ into $\mathfrak m$ (this is nontrivial). Since $\alpha(0)=x\in\mathfrak m$ and $G$ is a Lipschitz map by the uniqueness of the solutions of  ordinary differential equations we have $\alpha\subset\mathfrak m$. Hence ${\rm e}^{\alpha(1)}=qpq\in M$ and the claim follows.\end{proof} 

\begin{coro}\label{triple}
Assume $M={\rm exp}(\mathfrak m)\subset\Sigma$, and $\mathfrak m$ is as in the above theorem. Then $M$ is a closed convex submanifold.
\end{coro}
\begin{proof}
Take $p,q\in M$. Then $p={\rm e}^x$, $q={\rm e}^y$ with $x,y\in \mathfrak m$. If we put $r={\rm e}^{-x/2}{\rm e}^y{\rm e}^{-x/2}$,  then $r\in M$ because ${\rm e}^{-x/2}$ and ${\rm e}^y$ are in $M$. Moreover, $z=\ln(r)\in \mathfrak m$. But the unique geodesic of $\Sigma$ joining $p$ to $q$ is $\gamma(t)={\rm e}^{x/2}{\rm e}^{tz} {\rm e}^{x/2}$,  hence $\gamma\subset M$.
\end{proof}

\begin{coro}\label{algebra}
Assume ${\mathfrak m}\subset{\cal H}_{\mathbb R}$ is a closed, commutative associative Banach subalgebra of ${\cal H}_{\mathbb C}$. Then the manifold $M=\mbox{\rm exp}(\mathfrak m)\subset\Sigma$ is a closed, convex and flat Riemannian submanifold. Moreover, $M$ is an open subset of $\mathfrak m$ and an abelian Banach-Lie group.
\end{coro}
\begin{proof}
The first assertion follows from the fact that $\mathfrak m$ is a Lie triple system. Curvature is given by commutators, hence $M$ is flat. Since $\mathfrak m$ is a closed subalgebra, ${\rm e}^x=\sum \frac{x^n}{n!}\in \mathfrak m$ for any $x\in\mathfrak m$, so $M\subset \mathfrak m$. That $M$ is open follows from the fact that ${\rm exp}$ is a $C^{\omega}$ isomorphism (Corollary \ref{esiso}).\end{proof}

If $M$ is flat and geodesic at $p=1$, $T_1M=\mathfrak m$ is abelian (by Proposition \ref{abe}), therefore

\begin{coro}
Assume $M=\mbox{\rm exp}(\mathfrak m)$ is closed and flat. If $M$ is geodesic at $p=1$, then $M$ is a convex submanifold. Moreover, $M$ is an abelian Banach-Lie group and an open subset of $\mathfrak m$.
\end{coro}

We adopt the usual definition of a symmetric space \cite{sigurdur}:

\begin{defi}
A Riemann-Hilbert manifold $M$ is called a \emph{globally symmetric space} if each point $p\in M$ is an isolated fixed point of an involutive isometry $s_p:M\to M$. The map $s_p$ is called the \emph{geodesic symmetry}.
\end{defi}

\begin{teo} Assume $M={\rm exp}(\mathfrak m)$ is closed and convex. Then $M$ is a symmetric space; the geodesic symmetry at $p\in M$ is given by $s_p(q)=pq^{-1}p$ for any $q\in M$. In particular, $\Sigma$ is a symmetric space.
\end{teo}
\begin{proof}
Observe that, for $p={\rm e}^x$, $q={\rm e}^y$, $s_p(q)={\rm e}^x{\rm e}^{-y}{\rm e}^x$; this shows that $s_p$ maps $M$ into $M$. To prove that $s_p$ is an isometry, for any vector $v\in\mathfrak m$ consider the geodesic $\alpha_v$ of $M$ such that $\alpha(0)=q$ and $\dot\alpha(0)=v$. Then $\alpha(t)=q{\rm e}^{t\,q^{-1}v}$ and 
$$
d(s_p)_q(v)=\frac{d}{dt}|_{t=0} (s_p\circ \alpha_v)=-pq^{-1}vq^{-1}p.
$$
Since $M$ has the induced metric, $\|pq^{-1}vq^{-1}p\|^2_{pq^{-1}p}=\|v\|^2_q$  by  Lemma \ref{invar} (with $g=pq^{-1}$). In particular, $d_ps_p=-id$, so $p$ is an isolated fixed point of $s_p$ for any $p\in M$.
\end{proof}

\smallskip

Theorem \ref{pharpe} and its corollaries imply that $\Sigma$ (as any symmetric space) contains plenty of convex sets; in particular

\begin{rem} We can embed isometrically any $k$-dimensional plane in $\Sigma\;$ as a convex closed submanifold: take an orthonormal set of $k$ commuting operators (for instance, fix an orthonormal  basis  $\{e_i\}_{i\in\mathbb M}$ of $H$ and take $p_i=e_i\otimes e_i$, $i=1,\cdots, k$), and consider the exponential of the linear span of this set. In the language of symmetric spaces, we are saying that $rank\left(\Sigma\right)=+\infty$.\end{rem}

Let $I(M)$ be the group of isometries of a submanifold $M$.

\begin{teo}\label{isometrias}
If the submanifold $M={\rm exp}(\mathfrak m)$ is closed and convex, then $I(M)$ acts transitively on $M$.
\end{teo}
\begin{proof}
Take  $p={\rm e}^x$, $q={\rm e}^y$ two points in $M$, $v=p\ln(p^{-1}q)$ and $\gamma(t)=p{\rm e}^{t\,p^{-1}v}$ the geodesic joining $p$ to $q$. Note that $q=\gamma(1)=p{\rm e}^{\,p^{-1}v}={\rm e}^{\,vp^{-1}}p$. 
Consider the curve of isometries $\varphi_t=s_{\gamma(t/2)}\circ s_p$. Then
$$
\varphi_1(p)={\rm e}^{\frac12 v{\rm e}^{-x}}{\rm e}^x{\rm e}^{-x}\,{\rm e}^{\frac12 v{\rm e}^{-x}}{\rm e}^x={\rm e}^{\, v{\rm e}^{-x}}{\rm e}^x=q.
\qedhere
$$
\end{proof}

\begin{rem}
Assume $M\subset\Sigma$ is closed and convex, and $1\in M$. Let $I(M)$ be the group of isometries of $M$. Then, since any isometry $\varphi$ is uniquely determined by its value at $1\in M$ and its differential $d\varphi_1$, the set $I(M)$ can be naturally embedded in a Banach space: take $\varphi\in I(M)$ and consider 
$$
{\overline{\varphi}}\;(q)=\varphi(1)^{-\frac12}\; \varphi(q)\; \varphi(1)^{-\frac12}.
$$
Note that $d{\overline{\varphi}}_1$ is a unitary operator of  $T_1M=\mathfrak m$ (with the  natural Hilbert-space structure), so there is an inclusion $J:I(M)\hookrightarrow M\times {\cal U}({{\rm B}(\mathfrak m)})$ given by the map $\varphi\mapsto (\varphi(1),d{\overline\varphi}_1)$. On the other hand, for a given pair $(p,u)\in M\times {\cal U}({{\rm B}(\mathfrak m)})$, put $\varphi(\mbox{e}^x)= p^{\frac12} \exp(u(x)) p^{\frac12}$, $(x\in \mathfrak m)$. It is not hard to see that $\varphi$ is an isometry of $M$ which maps $1$ to $p$, such that $d{\overline\varphi}_1=u$ . Hence we may identify $I(M)\simeq M\times {\cal U}({{\rm B}(\mathfrak m)})$.
\end{rem}

\begin{rem}\label{tangente}
If $M={\rm exp}(\mathfrak m)$ is closed and convex, it is geodesic at any $p={\rm e}^x\in M$, so 
$$
T_pM={\rm Exp}_p^{-1}(M)= \{p^{\frac12}\ln(p^{-\frac12}\,q\,p^{-\frac12})p^{\frac12}: q\in M \}
$$
(see Remark \ref{exponencial}). Since $p^{\frac12}={\rm e}^{x/2}\in \mathfrak m$, using Theorem \ref{pharpe} we obtain the  identification $\displaystyle T_pM=p^{\frac12}\left(T_1M\right) p^{\frac12}=p^{\frac12}\, \mathfrak m\, p^{\frac12}$. It also follows easily that an operator $v\in {\cal H}_{\mathbb R}$ is orthogonal to $M$ at $p$ (that is, $v\in T_pM^{\perp}$) if and only if 
$$
<p^{-\frac12}\, z \, p^{-\frac12}\;,\; v>_{_2}=<p^{-\frac12}\, v \, p^{-\frac12},z>_{_2}=0\quad \mbox{ for any } z\in \mathfrak m.
$$
In particular, $\displaystyle T_1M^{\perp}=\mathfrak m^{\perp}=\{v\in {\cal H}_{\mathbb R}\;:\; <v,z>_{_2}=0\;\mbox{ for any }\;z\in \mathfrak m\}$. Note that, when $\mathfrak m$ is a closed commutative \emph{associative subalgebra} of operators,  $y\mapsto p^{\frac12} y p^{\frac 12}$ is a linear isomorphism of $\mathfrak m$;  in this case $\;T_pM=\mathfrak m=T_1M\;$ for any $p\in M$. This last assertion also follows easily from Corollary \ref{algebra}, and clearly $\displaystyle T_pM^{\perp}=T_1M^{\perp}=\mathfrak m^{\perp}$ in this case.
\end{rem}

\begin{rem}\label{transporte}
Assume $M\subset\Sigma$ is a convex submanifold. If the curve $\gamma$ is the geodesic joining $p$ to $q$, then  the isometry $\varphi_t=s_{\gamma(t/2)}\circ s_p$ translates along $\gamma$, namely
$$
\varphi_t(\gamma(s))=p\;{\rm e}^{\, \frac t2  p^{-1}v} \, p^{-1} \, p\; {\rm e}^{s p^{-1} v} \, p^{-1}\, p \;{\rm e}^{\, \frac t2 p^{-1}v}=\qquad\qquad\qquad\quad
$$
$$
\quad =p\;{\rm e}^{\, \frac t2  p^{-1}v} \, {\rm e}^{s p^{-1} v} \,  \;{\rm e}^{\, \frac t2 p^{-1}v}=p\;{\rm e}^{(s+t) p^{-1} v} =\gamma(s+t).
$$
In particular, $\varphi_1(p)=q$. Now take any tangent vector $w\in T_{\gamma(s)}M$, and let 
$$
w(t)=(d\varphi_t)_{\gamma(s)}(w)={\rm e}^{\frac t2 v p^{-1} }w\, {\rm e}^{\frac t2 p^{-1} v}.
$$
It follows from a straightforward computation using equation (\ref{covariant}) of Section \ref{intro} that $w(t)$ is the parallel translation of $w$ from $\gamma(s)$ to $\gamma(s+t)$; namely $\nabla_{\dot\gamma}\;w\equiv 0$. We conclude that the linear map $(d\varphi_t)_{\gamma(s)}:T_{\gamma(s)}M \to T_{\gamma(s+t)}M$ gives parallel translation along $\gamma$, \textit{i.e}  $(d\varphi_t)_{\gamma(s)}=P^{t+s}_s(\gamma)$. In particular, since $q=\gamma(1)=p^{\frac12}\;{\rm e}^{p^{-\frac12} vp^{\frac12}}p^{\frac 12}$, the map
$$
P^q_p : w\mapsto p^{\frac 12}(p^{-\frac12}qp^{-\frac12})^{\frac 12}p^{-\frac12} \; w \; p^{-\frac 12}(p^{-\frac12}qp^{-\frac12})^{\frac 12}p^{\frac12} 
$$
gives parallel translation from $T_pM$ to $T_qM$. See also Theorem \ref{homo}.\end{rem}

\subsubsection{Examples of convex sets}\label{examples}

\begin{enumerate}
\item For any subspace $\mathfrak s\subset  {\cal H}_{\mathbb R}$, ${\mathfrak m}_{\mathfrak s}=\{x\in {\cal H}_{\mathbb R}\;:\;[x,y]=0\;\forall \;y\in\mathfrak s\}$
 is a Lie triple system.
\item In particular, for any $y\in {\cal H}_{\mathbb R}$, ${\mathfrak m}_y=\{x\in {\cal H}_{\mathbb R}\;:\;[x,y]=0\}$ is a Lie triple system.
\item The family of operators in ${\cal H}_{\mathbb R}$ which act as endomorphisms of a closed subspace  $S\subset H$ form a Lie triple system in ${\cal H}_{\mathbb R}$.
\item Any norm closed commutative associative subalgebra of ${\cal H}_{\mathbb R}$, closed under the usual involution of operators, is a Lie triple system. In particular
\begin{enumerate}
\item The diagonal operators (see Section \ref{aplidiag}). This is a maximal abelian closed subspace of ${\cal H}_{\mathbb R}$, hence the manifold $\Delta$ (which is the exponential of this set) is a maximal flat submanifold of $\Sigma$.
\item The scalar manifold $\Lambda=\{\lambda 1:\lambda\in\mathbb R_{>0}\}\subset\Sigma$ is the exponential of the Lie triple system $\mathbb R\, 1\subset {\cal H}_{\mathbb R}$.
\item For fixed $a\in {{\sf HS}}^h$, the real part of the closed algebra generated by $a$, which is the closure in the 2-norm of the set of polynomials in $a$, is a Lie triple system.
\end{enumerate}
\item The real part of any Banach-Lie subalgebra of ${\cal H}_{\mathbb C}$ is a Lie triple system (in particular: the real part of any associative Banach subalgebra).
\end{enumerate}

\subsection{Convex manifolds as homogeneous manifolds}\label{homogeneo}

The results of this section are related to those of Sections 3 and 7 of Chapter IV in \cite{sigurdur}. See also Theorem 5.5 in \cite{neeb} for a proof of the existence of smooth  polar decompositons in the (broader) Banach-Finsler context.

\begin{defi}
Let ${\cal H}_{\mathbb C}^{\bullet}$ be the group of invertible elements in ${\cal H}_{\mathbb C}$. This group has a natural structure of manifold as an open set of the associative Banach algebra ${\cal H}_{\mathbb C}$; it is a Banach-Lie group with Banach-Lie algebra ${\cal H}_{\mathbb C}$.

Let ${\cal U}({\cal H}_{\mathbb C})$ stand for the unitary elements of the involutive Banach algebra ${\cal H}_{\mathbb C}$, namely the set of $u\in {\cal H}_{\mathbb C}^{\bullet}$ such that $u^*=u^{-1}$. It is a real Banach-Lie subgroup of ${\cal H}_{\mathbb C}^{\bullet}$ with Lie algebra $i{\cal H}_{\mathbb R}$.

Let $G$ be a connected abstract subgroup of ${\cal H}_{\mathbb C}^{\bullet}$. We say that $G$ is a self-adjoint subgroup of ${   {\cal H}_{\mathbb C}  }^{\bullet}$ if $g^*\in G$ whenever $g\in G$ (for short, $G^*=G$). Note that a connected Banach-Lie group $G$ is self-adjoint if and only if $\mathfrak g^*=\mathfrak g$, where $\mathfrak g$ denotes the Banach-Lie algebra of $G$. 

If $\mathfrak a\subset{\cal H}_{\mathbb C}$ is a linear space over $\mathbb R$, let $
[\mathfrak a,\mathfrak a]=\overline{span_{\mathbb R}\{[a,b]:a,b\in\mathfrak a\}}
$, where the bar denotes closure in the norm of the Banach algebra ${\cal H}_{\mathbb C}$.

If $A\subset {\cal H}_{\mathbb C}^{\bullet}$ is a set, $\left\langle A \right\rangle$ will denote the abstract subgroup generated by $A$ (the group whose elements are the inverses and the finite products of elements in $A$).

Let $\mid x\mid=(xx^*)^{\frac12}=\exp(\frac12 \ln(xx^*))$ for $x\in {\cal H}_{\mathbb C}$. Since ${\cal H}_{\mathbb C}$ is an involutive Banach algebra, $\mid x\mid\in \Sigma\subset{\cal H}_{\mathbb C}^{\bullet}$ if $x\in {\cal H}_{\mathbb C}$. 
\end{defi}

\begin{rem}\label{conexo}
The group ${   {\cal H}_{\mathbb C}  }^{\bullet}$, having the homotopy type of the inductive limit of the groups $GL(n,\mathbb C)$ (see \cite{pharpe}, Section II.6) is connected; moreover, there is a homotopy class equivalence
$$
{   {\cal H}_{\mathbb C}  }^{\bullet}\simeq S^1\times S^1 \times SU(\infty).
$$
Here $SU(\infty)$ stands for the inductive limit of the groups $SU(n,\mathbb C)$.
\end{rem}

\begin{prop}\label{integra}
Let $\mathfrak g\subset {\cal H}_{\mathbb C}$ be a closed real Banach-Lie subalgebra. Then $G=\left\langle \exp(\mathfrak g) \right\rangle$ admits a topology and a smooth structure such that $G$ is a connected real Banach-Lie group and $\mathfrak g=T_1G$ is the Banach-Lie algebra of $G$. The inclusion $G\hookrightarrow  {\cal H}_{\mathbb C}^{\bullet}$
 is a smooth inmersion and the exponential map of $G$ is given by the usual exponential of ${\cal H}_{\mathbb C}$. The topology on $G$ might be strictly finer than the topology of $ {\cal H}_{\mathbb C}^{\bullet} $.
\end{prop}
\begin{proof}
Since ${\cal H}_{\mathbb C}$ is a Hilbert space, the Banach-Lie subalgebra admits a suplement. By Theorem 5.4 of Chapter VI in \cite{lang}, there exists an integral manifold $H\stackrel{j}{\hookrightarrow}{\cal H}_{\mathbb C}^{\bullet}$ for the subbundle $\{g\mathfrak g\}_{g\in {\cal H}_{\mathbb C}^{\bullet}}$. The manifold $H$ is connected, and a Banach-Lie group with $dj_1(T_1H)=\mathfrak g$. Since $j$ is a smooth homomorphism of Banach-Lie groups, we have  $j\circ Exp^H=\exp\circ dj_1$. The other assertions follow from this identity because $G=\left\langle \exp(\mathfrak g) \right\rangle=\left\langle j\circ Exp^H(T_1H) \right\rangle=j(H)$.
\end{proof}

\begin{teo}\label{homo}
Let $G=\left\langle \exp(\mathfrak g) \right\rangle\subset {   {\cal H}_{\mathbb C}  }^{\bullet}$ be a connected self-adjoint Banach-Lie group with Banach-Lie algebra $\mathfrak g\subset{\cal H}_{\mathbb C}$. Let $P$ be the analytic map $g\mapsto gg^*$, $P:G\to G$. Let $\mathfrak k=\ker(dP_1)$, $\mathfrak m=\mathrm{Ran}(dP_1)$. Let $M_G=\exp(\mathfrak m)$,  $K=G\cap {\cal U}({\cal H}_{\mathbb C})=P^{-1}(1)$.  Then 
\begin{enumerate}
\item The set $\mathfrak m$ is a closed Lie triple system in ${\cal H}_{\mathbb R}$. We have $[\mathfrak m,\mathfrak m]\subset \mathfrak k$, $[\mathfrak m,\mathfrak k]\subset \mathfrak m$, $[\mathfrak k,\mathfrak k]\subset \mathfrak k$ and $\mathfrak g=\mathfrak m\oplus\mathfrak k$. In particular, $\mathfrak k$ is a Banach-Lie subalgebra of $\mathfrak g$ (and of $i{\cal H}_{\mathbb R}$ also).
\item $P(G)=M_G$, and $M_G$ is a geodesically convex submanifold of $\Sigma$.
\item For any $g=\mid g\mid u_g\in G$ (polar decomposition), we have $\mid g\mid\in M_G$ and $u_g\in K$.
\item Let $g\in G$, $p\in M_G$, $I_g(p)=gpg^*$. Then $I_g  \in I(M_G)$. If $g=p^{\frac 12}(p^{-\frac12}qp^{-\frac12})^{\frac 12}p^{-\frac12} \in G$, then $I_g(p)=q$, namely $G$ acts isometrically and transitively on $M_G$.
\item Let $u\in K$ and $x\in\mathfrak m$ (resp. ${\mathfrak m}^{\perp}$). Then $I_u(x)=uxu^*\in \mathfrak m$ (resp. ${\mathfrak m}^{\perp}$). If $p,q\in M_G$ then $I_{p}$ maps $T_qM_G$ (resp. ${T_qM_G}^{\perp}$) isometrically onto $T_{  I_p(q) }M_G$ (resp. ${T_{I_p(q)}M_G}^{\perp}$).
\item The group $K$ is a Banach-Lie subgroup of $G$ with Lie algebra $\mathfrak k$.
\item $G\simeq M_G\times K$ as Hilbert manifolds. In particular $K$ is connected and $G/K\simeq M_G$.
\end{enumerate}
\end{teo}
\begin{proof}
1. Note that $dP_1(x)=x+x^*$, hence $\mathfrak k=\{x\in\mathfrak g: x^*=-x \}$ which is certainly a closed Lie algebra. Note also that $\mathfrak m=\{x\in\mathfrak g: x^*=x \}$ is a Lie triple system; it is closed because $x\mapsto x^*$ is an isometric automorphism of ${\cal H}_{\mathbb C}$. Since $[x,y]=xy-yx$ is self-adjoint whenever $x$ is self-adjoint and $y$ is skew-adjoint, the other assertions are clear.

2. Cleary $P(G)\supseteq\exp(\mathfrak m)$ because ${\rm e}^x=P({\rm e}^{x/2})$. On the other hand, since $\mathfrak g$ splits, there exist neighbourhoods of zero $U_{\mathfrak m}\subset\mathfrak m$ and $U_{\mathfrak k}\subset \mathfrak k$ such that the map $x_{\mathfrak m}+y_{\mathfrak k}\mapsto {\rm e}^{x_{\mathfrak m}}{\rm e}^{y_{\mathfrak k}}$ is an isomorphism from $U_{\mathfrak m}\oplus U_{\mathfrak k}$ onto an open neighbourhood $V$ of $1\in G$. Then  $\left\langle V\right\rangle$ is open (and closed) in $G$ and so is all of $G$. Hence, for any $g\in G$, $g=({\rm e}^{x_1}{\rm e}^{y_1})^{\alpha_1}\cdots ({\rm e}^{x_n}{\rm e}^{y_n})^{\alpha_n}$ for self-adjoint $x_i\in U_{\mathfrak m}$, skew-adjoint $y_i\in U_{\mathfrak k}$, and $\alpha_i={\underline{+}1}$. Now ${\rm e}^x{\rm e}^y{\rm e}^x\in \exp(\mathfrak m)$ whenever $x,y\in \mathfrak m$ (Theorem \ref{pharpe}), and inspection of the expression for $P(g)=gg^*$ shows that $P(g)$ lies in $\exp(\mathfrak m)$ if  ${\rm e}^y{\rm e}^x{\rm e}^{-y}\in \exp(\mathfrak m)$ whenever $x\in\mathfrak m$ and $y\in \mathfrak k$. Equivalently, we have to show that $Ad( {\rm e}^y  )$ maps $\mathfrak m$ into $\mathfrak m$; since $Ad( {\rm e}^y )={\rm e}^{ad (y)}$, it suffices to show that $ad(y): x \mapsto [y,x]$ maps $\mathfrak m$ into $\mathfrak m$, and this follows from the previous assertion. The set $M_G=\exp(\mathfrak m)$ is a convex submanifold because $\mathfrak m$ is a closed Lie triple system (Corollary \ref{triple}). 

3. If $g\in G$, then $gg^*={\rm e}^{x_0}$ for some $x_0\in \mathfrak m$. This implies that $\mid g\mid ={\rm e}^{x_0/2}\in M_G\subset G$. Now we have $u_g=\mid g\mid^{-1} g\in G$, and clearly $u_g\in K$. 

4. If $p\in M_G$, then $p=P(g_o)=g_og_o^*$ for some $g_o\in G$. Then, if $g\in G$, $I_g(p)=gg_og_o^*g^*=P(gg_o)\in M_G$. Note that $I_g$  is an isometry of $M_G$, because $M_G$ has the induced metric, so Lemma \ref{invar} applies. 

5. If $x\in \mathfrak m$ and $u\in K$, then ${\rm e}^x\in M_G$ hence $u{\rm e}^xu^*=\exp(I_u(\mbox{e}^x))\in M_G$. Hence $uxu^*=\ln(u{\rm e}^xu^*)\in \mathfrak m$. Since $<I_u(y),x>_{_2}=<y,I_{u^*}(x)>_{_2}$ (see Remark \ref{darvuelta}), we obtain the proof of the assertion concerning ${\mathfrak m}^{\perp}$.

Clearly $I_{p}$ maps $T_qM_G$ isometrically onto $T_{I_p(q)}M_G$. Assume now $w\in {T_qM_G}^{\perp}=q^{\frac12}{\mathfrak m}^{\perp} q^{\frac12}$ (see Remark \ref{tangente}). If $u=(pqp)^{\frac12}p^{-1}q^{-\frac12}$, then $u\in G$ and $uu^*=1$, hence $w_0= u (q^{-\frac12} w q^{-\frac12}) u^*\in \mathfrak m^{\perp}$ by the previous assertion. Then
$I_{p}(w)=pwp= (pqp)^{\frac12} w_0 (pqp)^{\frac12} \in {T_{I_p(q)} M_G}^{\perp}$.

6. The previous items show that $P:G\to M_G$ is surjective. Now $dP_g:g.\mathfrak g\to T_{gg^*}M_G$ is given by $g.x\mapsto I_g(x+x^*)$. Clearly this map has split kernel $g\mathfrak k$. Let $g=\mid g\mid u_g$ as above. For $z\in T_{gg^*}M_G$ we have, by Remark \ref{tangente}, $z=(gg*)^{\frac12}w(gg*)^{\frac12}=I_{\mid g \mid}(w)$ for some $w\in\mathfrak m$. Let $x=I_{u_g^*}(w/2)$, then $x\in\mathfrak m\subset\mathfrak g$ and $dP_g(gx)=z$. Hence the group $K=P^{-1}(1)$ is a submanifold of $G$ because $P:G\to M_G$ is a submersion (Proposition 2.3 of Chapter II in \cite{lang}).

7. The map $T:M_G\times K\to G$ given by $T(p,u)=pu$ is clearly smooth and it is a bijection by the statements above. The inverse is given by $g\mapsto (\mid g\mid,\mid g\mid^{-1}g)$; since $\mid g\mid=\exp(\frac12 \ln(gg^*))$, the map $T$ is a diffeomorphism.
\end{proof}

\begin{rem}\label{geme}
For $M={\rm exp}(\mathfrak m)$ a convex closed manifold in $\Sigma$, consider $\mathfrak g_M=\mathfrak m\oplus [\mathfrak m,\mathfrak m]$. Then $\mathfrak g_M$ is a Banach-Lie subalgebra of ${\cal H}_{\mathbb C}$ due to the formal identity 
$$
[[x,y],[z,w]]+[z,[w,[x,y]]]+[w, [[x,y], z] ]=0
$$
and the fact that $\mathfrak m$ is a Lie triple system. Let $G_M=\langle \exp(\mathfrak g_M)\rangle$. Then $G_M$ is a connected Banach-Lie group with Banach-Lie algebra $\mathfrak g_M$ (Proposition \ref{integra}). Since $(a+[b,c])^*=a+[c,b]$ for any $a,b,c\in \mathfrak m$, then $M\subset G_M$ and $G_M^*=G_M$. It is also clear that $\mathfrak k=[\mathfrak m,\mathfrak m]$ ($\mathfrak k$ as in Theorem \ref{homo}). The elements of $M$ are indeed the positive elements of $G_M$, and the elements of the stabilizer of $1$ are the unitary operators of $G_M$.  Note that $G_M$ is a submanifold of ${\cal H}_{\mathbb C}^{\bullet}$ if and only if $K$ is a submanifold of ${\cal U}({\cal H}_{\mathbb C})$.
\end{rem}

When $\mathfrak m$ is a commutative associative subalgebra, we have $\mathfrak g_M=\mathfrak m$ and also $G_M=M\subset \mathfrak m$ is an open set (in particular $G_M$ is a submanifold of ${\cal H}_{\mathbb C}^{\bullet}$). In any case $\mathfrak m=Z(\mathfrak m)\oplus Z(\mathfrak m)^{\perp}=\mathfrak m_0\oplus \mathfrak m_1$ (here $Z(\mathfrak m)$ denotes the set $\{x\in\mathfrak m: [x,y]=0\quad\forall y\in\mathfrak m\}$), and $M=\exp(\mathfrak m_0)\exp(\mathfrak m_1)\simeq M_0\times M_1$ where $M_i$ are convex and closed, hence $G_M\simeq M_0\times G_{M_1}$. Since $<x,[a,[b,c]]>_{_2}=<c,[b,[a,x]]>_{_2}$ for any $a,b,c,x\in\mathfrak m$, it is easy to see that $Z(\mathfrak m)
=[\mathfrak m,[\mathfrak m,\mathfrak m]]^{\perp}$.

\medskip

The results above assert that, for a given convex submanifold $M=\exp(\mathfrak m)$, we have $M_{G_M}=M$. On the other hand, for a given connected involutive Banach-Lie subgroup $G$, we have $G_{M_G}\subset G$, though in general $[\mathfrak m,\mathfrak m]$ can be strictly smaller than $\mathfrak k$, so the other inclusion does not necessarily hold. The equality holds iff $\mathfrak k$ is semi-simple, \textit{i.e.} $[\mathfrak k,\mathfrak k]=\mathfrak k$ (equivalently, if $Z(\mathfrak k)=0$).

\medskip

It is a well known result (see \cite{pharpe}, p.42) that $\overline{[{{\sf HS}},{{\sf HS}}]}={{\sf HS}}$ and $\overline{[{{\sf HS}}^h,{{\sf HS}}^h]}=i{{\sf HS}}^h$. Therefore taking $\mathfrak m={{\sf HS}}^h$, we get $\mathfrak k=i{{\sf HS}}^h$, and then $\mathfrak g_M={{\sf HS}}$. This implies $
G_{\Sigma/{\mathbb R}}= {   {\cal H}_{\mathbb C}  }^{\bullet}/{\mathbb C^{\times}1}$. Clearly $P(G_{\Sigma}) = P({   {\cal H}_{\mathbb C}  }^{\bullet})=\Sigma$, because any positive invertible operator has an invertible square root. On the other hand it is clear that the isotropy group $K$ equals ${\cal U}\left({{\cal H}_{\mathbb C}}\right)$ (the unitary group of ${\cal H}_{\mathbb C}$). So there is an analytic isomorphism given by polar decomposition: $\Sigma\simeq {   {\cal H}_{\mathbb C}  }^{\bullet}/{\cal U}\left({{\cal H}_{\mathbb C}}\right)$. The manifold of positive invertible operators $\Sigma$ is an homogeneous manifold for the group of invertible operators ${   {\cal H}_{\mathbb C}  }^{\bullet}$, which acts isometrically and transitively on $\Sigma$. This last statement is well known, and Theorem \ref{homo} can be read as a natural generalization.

\section{Projecting to closed convex submanifolds}\label{elteo}

We refer the reader to \cite{lang} for the first and second variation formulas.

\begin{prop}\label{triang} Let $M$ be a convex subset of $\;\Sigma$, and let $p\in \Sigma$. Then there is at most one normal geodesic $\gamma$ of $\Sigma$ joining $p$ and $M$ such that $L(\gamma)=\mbox{\rm dist}\left(p,{M}\right)$. In other words, there is at most one point $q\in{M}$ such that  $\mbox{\rm dist}(p,q)=\mbox{\rm dist}\left(p,{M}\right)$.
\end{prop}
\begin{proof} Suppose there are two such points, $q$ and $r\in M$, joined by a geodesic $\gamma_3\in M$, such that $L(\gamma_1)=\mbox{\rm dist}(p,q)=L(\gamma_2)=\mbox{\rm dist}(p,r)=d\left(p,{M}\right)$. We construct a proper variation of $\gamma\equiv\gamma_1$, which we call $\Gamma_s$. The construction follows the figure below,  where $\sigma(s,t):=\sigma_s (t)$ is the geodesic joining $p$ with $\gamma_3(s)$.
\\
\begin{minipage}{7cm}
\hspace*{.2cm}\centerline{\includegraphics[totalheight=4cm,width=7cm]{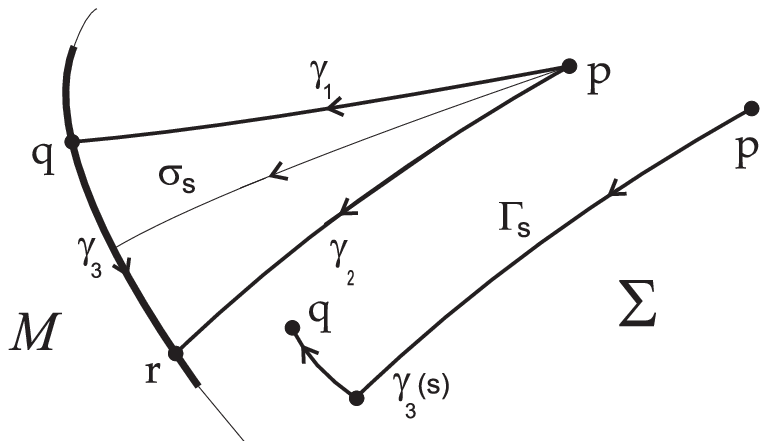}}
\end{minipage}

\hspace*{8.1cm}
\begin{minipage}{6cm}
\vspace*{-5cm}
$$
\Gamma(s,t)=    
 \left\{
 \begin{array}{lll}
 \sigma(s,t) &  0\le t\le 1 \\
 \gamma_3\left( s(2-t) \right) & 1 \le t \le 2\\
 \end{array}
 \right. ,
$$
\end{minipage}
\vspace*{-.4cm}
$$
\gamma(t)=\Gamma(0,t)=    
  \left\{
\begin{array}{lll}
\gamma_1(t) & \textnormal{if}& 0\le t\le 1 \\
q &\textnormal{if}&1 \le t \le 2\\
\end{array}
\right.,\mbox{ so that }\;  {\dot \gamma  (t)}=
\left\{
\begin{array}{lll}
{\dot \gamma_1  (t)} & \textnormal{if}& 0\le t\le 1 \\
0 &\textnormal{if}&1 \le t \le 2\\
\end{array}
\right. .
$$
Also note that the variation vector field (which is a Jacobi field for $\gamma$) is given by
$$
V(t)=\frac{\partial\Gamma}{\partial s}(t,0)=
\left\{
\begin{array}{ll}
\frac{\partial\sigma}{\partial s}(t,0) &  0\le t\le 1 \\
(2-t){\dot \gamma_3(0)} &1 \le t \le 2\\
\end{array}
\right. .
$$
If $\Delta_i {\dot \gamma}$ denotes the jump of the tangent vector field to $\gamma$ at $t_i$, namely ${\dot \gamma}(t_i^+)-{\dot \gamma}(t_i^-)$, and $\Gamma$ is a proper variation of $\gamma$, then the first variation formula for the curve $\gamma:[0,2]\to \Sigma$ reads
$$
\|\dot\gamma \|_{_{\gamma}} \;\frac{d}{ds} | _{_{s=0^+}} L\left(\Gamma_s\right) = -\int_0^2 <V(t),D_t{\dot \gamma(t)}>_{_{ \gamma(t) }}dt\; - \sum_{i=1}^{k-1}   <V(t_i),\Delta_i {\dot \gamma}>_{_{ \gamma(t_i) }} .
$$
In this case, $D_t{\dot \gamma}$ is zero in the whole interval $[0,2]$, because $\gamma$ consists  (piecewise) of  geodesics. The jump points are $t_0=0$, $t_1=1$ and $t_2=2$, so the formula reduces to
$$
<{\dot \gamma_3(0)},{\dot \gamma_1(1)}>_{_q}= \frac{d}{ds} | _{_{s=0^+}} L\left(\Gamma_s\right) \|\dot\gamma \|_{_{\gamma}}.
$$
Recall that $\gamma_3\subset M$, and that $\gamma_1$ is minimizing. Then the right hand term is nonnegative, which proves that the angle between $\gamma_1$ and $\gamma_3$ at $q$ is bigger that $\pi/2$. With a similar argument, we deduce that the same holds for the angle between $\gamma_2$ and $\gamma_3$ at $r$. Hence, the sum of the three inner angles of this geodesic triangle is \emph{at least} $\pi$. Since the sum cannot exceed $\pi$ (see Corollary \ref{triangulines}), it follows that the angle subtended at $p$ must be zero, which proves that $\gamma_1$ and $\gamma_2$ are the same geodesic, and uniqueness follows. 
\end{proof}

\medskip

Next we consider the problem of existence of the minimizing geodesic.

\begin{prop} Let $M$ be a convex submanifold of $\Sigma$, and 
$p$ a point of $\Sigma$ not in $M$. Then the existence of a geodesic $\beta$ joining $p$ with ${M}$ such that $L(\beta)=\mbox{\rm dist}(p,M)$ is equivalent to the existence of a geodesic $\gamma$ joining $p$ with ${M}$ with the property that $\gamma$ is orthogonal to ${M}$.\end{prop}
\begin{proof} In fact, the existence of such a geodesic $\beta$ is equivalent to the existence of a point $q\in {M}$ such that $\mbox{\rm dist}(p,{M})=\mbox{\rm dist}(p,q)$. We will show that  if $q\in{M}$ is a point such that $\gamma_{qp}$ is orthogonal to ${M}$ at $q$, then $\mbox{\rm dist}(q,p)=\mbox{\rm dist}(M,p)$. The other implication follows from the uniqueness theorem above. Consider the geodesic triangle generated by $p, q$ and $d$, where $d$ is any point in ${M}$ different from $q$. 
Since $\gamma_{qp}$ is orthogonal to $T_q{M}$, it is orthogonal to $\gamma_{qd}$. Then, by virtue of the Hyperbolic Cosine Law (equation (\ref{hypcos}) in Section \ref{desiguals}), ${\rm dist}(d,p)^2=L(\gamma_{dp})^2\ge L(\gamma_{qp})^2+L(\gamma_{qd})^2>L(\gamma_{qp})^2={\rm dist}(q,p)^2.$\end{proof}

This last proposition raises the following question: is the normal bundle $NM$ of ${M}$ diffeomorphic to $\Sigma$, \emph{via} the exponential map?

\begin{lema}\label{tubo}
Let $M$ be a convex and closed submanifold. Let $E:N{M}\to \Sigma$ be the map $(q,v)\mapsto \mbox{\rm Exp}_q(v)$. For $\varepsilon>0$, put $NM_{\varepsilon}=\{(p,v)\in NM:\|v\|_{_p}<\epsilon\}$ and  $\Omega_{\varepsilon}=E(NM_{\varepsilon})$. Then $E$ is injective and there exists $\varepsilon>0$ such that $E: NM_{\varepsilon}\to \Omega_{\varepsilon}$ is a $C^{\omega}$-diffeomorphism. The set $\Omega_{\varepsilon}$ is an open neighbourhood of $M$ in $\Sigma$.
\end{lema}
\begin{proof}
Let us prove first that $E$ is injective. Assume there exist $p,q\in M$, $v\in T_pM^{\perp}$, $w\in T_qM^{\perp}$ with ${\rm Exp}_p(v)={\rm Exp}_q(w)$. Naming $r$ to this point, consider the geodesic triangle in $\Sigma$ spanned by $p,q\in M$ and $r\in \Sigma$. The geodesic which joins $p$ to $r$ is clearly $\gamma_1(t)={\rm Exp}_p(tv)$, which is orthogonal to $M$ at $p$, and the same is true for $\gamma_2(t)={\rm Exp}_q(tw)$, which joins $q$ to $r$. Hence $p=q$ and $v=w$ because of Corollary \ref{triangulines}.

We may assume that $1\in M$. Since $E(q,v)=q{\rm e}^{q^{-1}v}$, the differential of $E$ at $(1,0)\in NM$ is the identity map because $T_1{M}\oplus  T_1{M}^{\perp}= T_1\Sigma$ and $d\exp_0=id$. The inverse mapping theorem (\cite{lang}, Theorem 5.2 of Chapter I) gives $C^{\omega}$-diffeomorphic neighbourhoods $U_0=\{(q,v)\in NM: {\rm dist}(q,1)<\varepsilon,\; \|v\|_{_q}<\varepsilon\}\subset NM_{\varepsilon}$  
and $\Omega_0=E(U_0)\subset \Sigma$ respectively. For given $(p,v)\in NM_{\varepsilon}$, consider the isometry of $M$ given by  $\tilde{I_p}: x\mapsto p^{\frac12}x p^{\frac12}$, and note that $\tilde{I_p}(1)=p$. If $(q,w)\in U_0$, then clearly $\tilde{I_p}(q)\in M$. Moreover, $\tilde{I_p}(w)\in T_{\tilde{I_p}(q)}M^{\perp}$ by Theorem \ref{homo}, hence $U_p=(\tilde{I_p}\times \tilde{I_p})(U_0)$ is an open neighbourhood of $(p,v)$  in $NM_{\varepsilon}$ diffeomorphic to $U_0$. Now $E\mid_{U_p}:U_p\to E(U_p)$ is a diffeomorphism, because a straightforward computation shows that $E\mid_{U_p}=   \tilde{I_p}\circ E\circ(\tilde{I_p}\times \tilde{I_p})^{-1}$.
\end{proof}

\begin{rem}
Clearly $E(NM)\subset\Sigma$ is the set of points $p\in\Sigma$ with the following property: there is a point $q\in {M}$ such that $\mbox{\rm dist}(q,p)=\mbox{\rm dist}\left({M},p\right)$. Note that the map $\Pi_{M}:  E(NM)\to {M}$, which assigns to  $p\in E(NM)$ the unique point $q\in {M}$ such that $\mbox{\rm dist}(q,p)=\mbox{\rm dist}\left({M},p\right)$, is surjective. This map is obtained \emph{via} a geodesic that joins $p$ and $M$, and this geodesic is orthogonal to $M$, therefore we call $\Pi_{M}(p)$ the \emph{foot of the perpendicular} from $M$ to $p$.
\end{rem}

\begin{lema}\label{convexa}
Let $p,q\in E(NM)$, and $\Pi_{M}(p)\ne \Pi_{M}(q)$. If $\gamma_p$ is a geodesic that joins $\Pi_{M}(p)$ to $p$ and $\gamma_q$ is a geodesic that joins $\Pi_{M}(q)$ to $q$, put $f(t)={\rm dist}\left(\gamma_p(t),\gamma_q(t)\right)$. Then the map $f:\mathbb R_{\ge 0}\to\mathbb R_{\ge 0}$ is increasing.
\end{lema}
\begin{proof} 
Since $f$ is a convex function (Corollary \ref{con}), it suffices to show that $f\; '(0^+)\ge 0$.

\hspace*{-.6cm}
\begin{minipage}{6.5cm}
\vspace*{.2cm}
\begin{spacing}{1.2}
Take a variation $\sigma(t,s)$, where $\sigma_t(s)$ is the geodesic joining $\gamma_p(t)$ to  $\gamma_q(t)$. Then 
$\sigma(t,0)=\gamma_p(t)$,
$\sigma(t,1)=\gamma_q(t)$, 
 and $\sigma(0,s)=\gamma(s)$ is  the geodesic joining $\Pi_{M}(p)$ to $\Pi_{M}(q)$ (which is contained in $M$ by virtue of the convexity). Note also that  $\sigma(1,s)$ is the geodesic joining $p$ to $q$. This construction is shown in the figure on the right.
\end{spacing}
\end{minipage}

\vspace*{-4.9cm}
\begin{minipage}{7cm}
\hspace*{7cm}\includegraphics[height=4.3875cm,width=5.805cm]{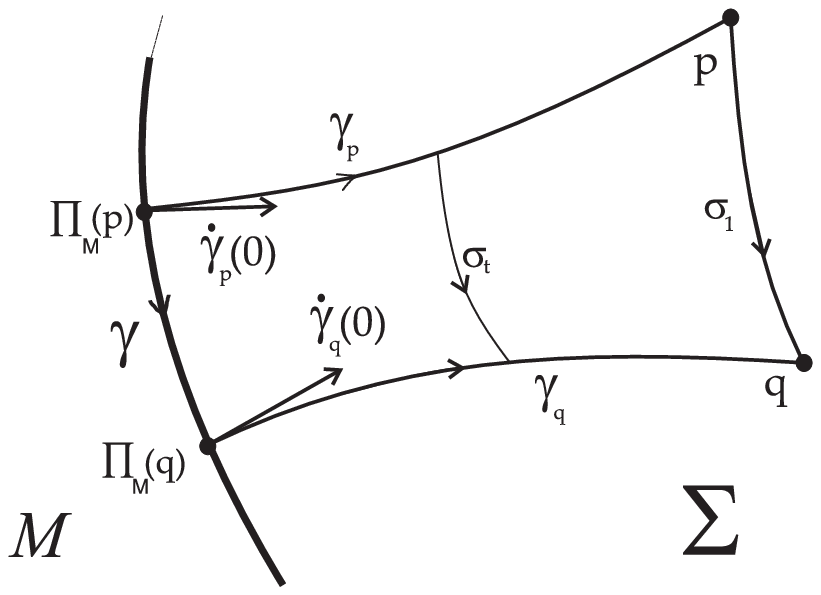} 
\end{minipage}

\vspace*{.2cm}

Note that $f(t)=L(\sigma_t)$. Put $V=\frac{d}{dt}|_{t=0}\sigma$. We apply the first variation formula to obtain
$$
\|\dot\gamma \|_{_{\gamma}} \;\frac{d}{dt} | _{_{t=0^+}} L\left(\sigma_t\right) = -\int_0^1 <V(s),D_s{\dot \gamma(s)}>_{_{\gamma(s)}}ds\;+\qquad\qquad\qquad\qquad
$$
$$
\qquad\qquad\qquad\qquad +<V(1),\dot\gamma(1)>_{_{\Pi_{M}(p)}}-<V(0),\dot\gamma(0)>_{\Pi_{M}(q)}.
$$
The fact that $\gamma$ is a geodesic reduces the formula to
$$
\|\dot\gamma \|_{_{\gamma}} \; f\;'(0^+)= -<V(1),-\dot\gamma(1)>_{_{\Pi_{M}(p)}}+<-V(0),\dot\gamma(0)>_{_{\Pi_{M}(q)}}.
$$
Note also that $V(0)={\dot\gamma_p}(0), $ $V(1)={\dot\gamma_q}(0) $. Recalling that the angles at $M$ are right angles, we obtain $f\; '(0^+)=0$. \end{proof}

\begin{teo}\label{contrae}The map $\,\Pi_{M}$ is a contraction, namely
$\mbox{\rm dist}\left(\Pi_{M}(p),\Pi_{M}(q)\right) \le \mbox{\rm dist}(p,q)
$.
\end{teo}
\begin{proof} 
We may assume again that $p,q\notin {M}$, and that $\Pi_{M}(p)\ne \Pi_{M}(q)$. In the notation of the lemma above, note that $f(0)=d\left(\Pi_{M}(p),\Pi_{M}(q)\right)$ and $f(1)=\mbox{\rm dist}(p,q)$; since $f$ is increasing, the assertion is proved.
\end{proof}

We want to prove that $E(NM)=\Sigma$. We will do this by proving that it is both open and closed in $\Sigma$. The following argument is similar to the one used by H. Porta and L. Recht in \cite{pr4}. 

\begin{lema}
For $\lambda\in [1,+\infty)$, put $\eta_{\lambda}:E(NM) \to E(NM)$, $\eta_{\lambda}({\rm Exp}_p(v))={\rm Exp}_p(\lambda v)$. Let $\Omega_{\varepsilon}$ be as in Lemma \ref{tubo}. Then $E(NM)=\displaystyle\mathop{\cup}_{\lambda\ge 1}\eta_{\lambda}(\Omega_{\epsilon})$, and each $\eta_{\lambda}:\Omega_{\epsilon}\to \Sigma$ is a $C^{\omega}$ diffeomorphism onto its open image.
\end{lema}

\begin{proof}
Clearly $\displaystyle\cup_{\lambda\ge 1}\eta_{\lambda}(\Omega_{\epsilon})\subset E(NM)$. Let us prove the other inclusion. First, if $r={\rm Exp}_p(v)$ with $\|v\|_{_p}<\epsilon$ then $r\in\Omega_{\epsilon}=\eta_1(\Omega_{\epsilon})$. Let us consider the case where $\|v\|_{_p}\ge \epsilon$; then $r={\rm Exp}_p(v)={\rm Exp}_p(\lambda w)$ with $\lambda=\frac{2\|v\|_{_p}}{\varepsilon}$ and $w=\frac{\varepsilon}{2\|v\|_{_p}} v$, so $r\in \eta_{\lambda}(\Omega_{\varepsilon})$ because $\|w\|_{_p}=\varepsilon/2<\varepsilon$ and $\lambda\ge 1$.

Assume that there exist $r_1,r_2\in \Omega_{\varepsilon}$ and $\lambda\ge 1$ such that $\eta_{\lambda}(r_1)=\eta_{\lambda}(r_2)$. That is, assume there exist $p,q\in M$, $v\in T_pM^{\perp}$, $w\in T_qM^{\perp}$ with $\|v\|_{_p}<\varepsilon$, $\|w\|_{_q}<\varepsilon$ and ${\rm Exp}_p(\lambda v)={\rm Exp}_q(\lambda w)$, namely $E(p,\lambda v)=E(q,\lambda w)$. Since $E$ is injective by Lemma \ref{tubo}, we have $p=q$ and $v=w$. This argument proves that the maps $\eta_{\lambda}$ are injective.

Next we show that, for any $\lambda\ge 1$ and $r\in \Omega_{\epsilon}$, $d\left(\eta_{\lambda}\right)_r:T_r\Sigma\to T_{\eta_{\lambda}(r)}\Sigma$ is a linear isomorphism, and this will prove the final assertion. Take $\alpha\subset \Omega_{\epsilon}$ a geodesic such that $\alpha(0)=r$ and $\dot\alpha(0)=x$. Since $\alpha$ is a geodesic, we have that ${\rm dist}(\alpha(t),r)=t\|\dot\alpha(0)\|_{_r}$ for $t\ge 0$ (see Section \ref{geods}).  Put $\beta=\eta_{\lambda}\circ\alpha$. Then $\beta(0)=\eta_{\lambda}(r)$ and $\dot\beta(0)=d\left(\eta_{\lambda}\right)_r (x)$. Clearly ${\rm dist}(\beta(t),\eta_{\lambda}(r)) \le L_0^t(\beta)=\int_0^t\|\dot\beta(s)\|_{_{\beta(s)}}ds$. On the other hand, $
{\rm dist}(\eta_{\lambda}(\alpha(t)), \eta_{\lambda}(r)) \ge {\rm dist}(\alpha(t),r)=t\|x\|_{_r}$ where the inequality is due to Lemma \ref{convexa}, because $\lambda\ge 1$. If we put together these two inequalities and divide by $t$, we get
$$
\frac1t \int_0^t\|\dot\beta(s)\|_{_{\beta(s)}}ds \ge \|x\|_{_r}.
$$
Taking limit for $t\to 0^+$ gives $ \| d\left(\eta_{\lambda}\right)_r (x)\|_{_{\eta_{\lambda}(r)}} \ge \|x\|_{_r}
$. Now put $A_{\lambda}=\tilde{I}_{\eta_{\lambda}(r)}^{-1}\circ d\left(\eta_{\lambda}\right)_r\circ \tilde{I}_r$, where $\tilde{I}_p: v\mapsto p^{\frac12}vp^{\frac12} $ are linear isomorphisms (see Lemma \ref{invar}). If we consider $A_{\lambda}:T_1\Sigma\to T_1\Sigma={\cal H}_{\mathbb R}$, the last inequality says that $\|A_{\lambda}(x)\|_{_2}\ge \|x\|_{_2}$ for any $x\in {\cal H}_{\mathbb R}$.

Clearly $\eta_1=id_{\Omega_{\epsilon}}$ and $d\left(\eta_{1}\right)_r=id_{T_r\Sigma}$. Since the map $(\lambda,r)\mapsto \eta_{\lambda}(r)$ is analytic from $\mathbb R_{>0}\times \Omega_{\epsilon}$ to $\Sigma$, there is an open neighbourhood of $1\in\mathbb R$ such that  $A_{\lambda}$ is an isomorphism. Assume $A_{\lambda}$ is invertible for $\lambda\in [1,m)$: then $\|A_{\lambda}^{-1}\|_{{\rm B}({\cal H}_{\mathbb R})}\le 1$ for any $\lambda \in [1,m)$. Since $A_m=\lim\limits_{\lambda \to m^-}A_{\lambda}$ (in the operator norm of ${\rm B}({\cal H}_{\mathbb R})$) and $\|A_mA_{\lambda}^{-1}-1\|\le \| A_m-A_{\lambda}\|<1$ if $\lambda$ is close enough to $m$, it follows that $A_mA_{\lambda}^{-1}$ is invertible, thus $A_m$ is invertible. Since the maps $\tilde{I}_p$ are isomorphisms, $d\left(\eta_{\lambda}\right)_r$ is an isomorphism for any $\lambda\ge 1$, and any $r\in \Omega_{\epsilon}$.\end{proof}

\begin{coro}
The set $E(NM)$ is open in $\Sigma$.
\end{coro}

\begin{teo}\label{main1} Let ${M}$ be a convex closed submanifold of $\Sigma$. Then for every point $p\in\Sigma$, there is a unique normal geodesic $\gamma_p$ joining $p$ to ${M}$ such that $L(\gamma_p)=\mbox{\rm dist}\left(p,{M}\right)$. This geodesic is orthogonal to ${M}$, and if $\Pi_{M}:\Sigma\to {M}$ is the map that assigns to $p$ the end-point of $\gamma_p$, then $\Pi_{M}$ is a contraction for the geodesic distance.\end{teo}
\begin{proof} The theorem will follow if we prove that $E(NM)=\Sigma$. Since $\Sigma$ is connected and $E(NM)$ is open, it suffices to prove that $E(NM)$ is also closed. Let $p\in \overline{E(NM)}$. There exist points $q_n\in {M}$, $v_n\in T_{q_n}M^{\perp}$ such that $\displaystyle p=\lim\limits_n p_n=\lim\limits_n \mbox{\rm Exp}_{q_n}(v_n)
$. Observe that $q_n=\Pi_{M}(p_n)$, so $\mbox{\rm dist}(q_n,q_m)\le \mbox{\rm dist}(p_n,p_m)$. Since $\{p_n\}$ converges to $p$, it is a Cauchy sequence. It follows that $\{q_n\}$ is also a Cauchy sequence. Since ${M}$ is closed (and therefore complete), there exists $q\in {M}$ such that $q=\lim\limits_n q_n$. We assert that $\mbox{\rm dist}(p,q)=\mbox{\rm dist}(p,{M})$. First note that $
\mbox{\rm dist}(p,q_n)\le \mbox{\rm dist}(p,p_n)+\mbox{\rm dist}(p_n,q_n)$ 
and $\mbox{\rm dist}(p_n,q_n)=\mbox{\rm dist}(p_n,{M})$, so $\displaystyle \mbox{\rm dist}(p,q_n)\le \mbox{\rm dist}(p,p_n)+\mbox{\rm dist}(p_n,{M})$. 
Taking limits gives $\mbox{\rm dist}(p,q)\le \mbox{\rm dist}(p,{M})$.\end{proof}

\smallskip

Note that $\Sigma$ decomposes as a direct product: with the contraction $\Pi_{M}$, we can decompose $\Sigma$ by picking, for fixed $p$, 
\begin{enumerate}
\item the unique point $q=\Pi_{M}(p)$ such that $\mbox{\rm dist}(p,q)= \mbox{\rm dist}(p,{M})$
\item a vector $v_p$ normal to $T_q{M}$ such that the geodesic in $\Sigma$ with initial velocity $v_p$ starting at $q$ passes through $p$; 
note that $v_p= \mbox{\rm Exp}^{-1}_{\Pi_ {{M}}(p)}(p),\; $ and also $\,\|v_p\|_{_q}=\mbox{\rm dist}(p,{M})$.
\end{enumerate}

\medskip

Since the exponential map is analytic on both of its variables, we get

\begin{teo} The map $p\mapsto \left(\Pi_{M}(p), v_p\right)$ is the inverse of the  map $(q, v_q)\mapsto \mbox{\rm Exp}_q(v_q),$ and gives a real-analytic isomorphism between the manifolds $\Sigma$ and $NM$.\end{teo}

\begin{teo}\label{gral1} Fix a closed convex submanifold ${M}$ of $\Sigma$. Let $a\in \Sigma$. Then there exist unique operators $c\in \Sigma$, $v\in {\cal H}_{\mathbb R}$ such that $c\in M$, $v\in T_c{M}^{\perp}$, and  $\displaystyle a=c\,\mbox{\rm e}^{c^{-1}v}$.
\end{teo}

Using the tools of Section \ref{themain}, we can write the factorization theorem in terms of intrinsic operator equations (see \cite{mostow} for the finite dimensional analogue):

\begin{teo}\label{gral2} Assume $\mathfrak m\subset {\cal H}_{\mathbb R}$ is a  Lie triple system. Then for any operator $a\in {\cal H}_{\mathbb R}$, there exist unique operators $x\in \mathfrak m$ and $v\in \mathfrak m^{\perp}$ such that the following decomposition holds:
$\displaystyle{\rm e}^a={\rm e}^{x}\,{\rm e}^v\,{\rm e}^{x}$. The map $d_a:{\cal H}_{\mathbb R}\to\mathbb R$, $d_a(y)= \|\ln({\rm e}^{a/2}{\rm e}^{-y}{\rm e}^{a/2})\|_{_2}$ has the operator $2x$ as its unique minimizer in $\mathfrak m$.
\end{teo}

\medskip

As a corollary, we obtain a polar decomposition relative to a convex submanifold.

\begin{teo}\label{mpolar}
Assume $M={\rm exp}(\mathfrak m)\subset\Sigma$ is a closed convex submanifold. Then for any $g\in {   {\cal H}_{\mathbb C}  }^{\bullet}$ there is a unique factorization of the form $g={\rm e}^x {\rm e}^v u$ where $x\in \mathfrak m$, $v\in \mathfrak m^{\perp}$ and $u\in {\cal U}({\cal H}_{\mathbb C})$ is a unitary operator. The map $g\mapsto ({\rm e}^x,{\rm e}^v,u)$ is an analytic bijection which gives an isomorphism
$$
{   {\cal H}_{\mathbb C}  }^{\bullet}\simeq M\times {\rm exp}(\mathfrak m^{\perp})\times {\cal U}({\cal H}_{\mathbb C}).
$$
\end{teo}
\begin{proof}
Since $gg^*\in\Sigma$, we can write $gg^*={\rm e}^{x}{\rm e}^{2v}{\rm e}^{x}$ with $x\in\mathfrak m$ and $v\in \mathfrak m^{\perp}$. If $u=({\rm e}^x{\rm e}^v)^{-1}g={\rm e}^{-v}{\rm e}^{-x}g$ we have $uu^*={\rm e}^{-v}{\rm e}^{-x}gg^*{\rm e}^{-x}{\rm e}^{-v}=1$ and also $u^*u=g^*{\rm e}^{-x}{\rm e}^{-v}{\rm e}^{-v}{\rm e}^{-x}g=1$. Hence $u$ is a unitary operator and $g={\rm e}^x{\rm e}^vu$. This factorization is unique because if $g={\rm e}^{x_1}{\rm e}^{v_1}u_1={\rm e}^{x_2}{\rm e}^{v_2}u_2$, then $gg^*={\rm e}^{x_1}{\rm e}^{2v_1}{\rm e}^{x_1}={\rm e}^{x_2}{\rm e}^{2v_2}{\rm e}^{x_2}$, so $x_1=x_2$, $v_1=v_2$ and then $u_1=u_2$.
\end{proof}

\section{Projecting to the manifold of diagonal operators}\label{aplidiag}

\begin{lema}\label{ef} Let $\alpha,\beta\in\mathbb R$ and $a,b\in {\sf HS}^h$. Then
$$
\mbox{\rm Exp}_{\alpha+a}(\beta+b)=\alpha \;\mbox{\rm e}^{\beta/\alpha}+ k
$$
where $k$ is a self-adjoint Hilbert-Schmidt operator.\end{lema}
\begin{proof}
It is a straightforward computation:
$$
(\alpha+a) {\rm e}^{(\alpha+a)^{-1}(\beta+b)}=(\alpha+a)[1+ (\alpha +a)^{-1}(\beta+b)+\cdots ]
$$
$$
\qquad\quad\qquad\quad\qquad\qquad\; = (\alpha+a)[1+\beta/\alpha+\frac{1}{2} \left(\beta/\alpha\right)^2+\cdots +k].\qedhere
$$
\end{proof}

We need some remarks before we proceed. Fix an orthonormal basis $\{e_i\}_{i\in \mathbb N}$ of $H$.

\begin{enumerate}
\item Consider the diagonal manifold ${\Delta}\subset\Sigma$: 
$$
\Delta=\{ d+\alpha>0 : \alpha\in\mathbb R,\; d\mbox{ is a diagonal Hilbert-Schmidt operator}\}.
$$
It is closed and geodesically convex. This is due to the fact that the diagonal operators form a closed commutative associative subalgebra.
\item If $d_0\in\Delta$, then $\displaystyle T_{d_0}{\Delta}=\{\alpha+d;\;\alpha\in\mathbb R,\; d\in{\sf HS}\textnormal{ is diagonal and real} \}=T_1\Delta$ (see Remark \ref{tangente}).
\item Consider the map $A\mapsto A^D$ = the diagonal part of $A$. Then
\begin{enumerate}
\item For Hilbert-Schmidt operators we have $A^D=\sum_i p_i A p_i$ where convergence is in the 2-norm (and hence in the operator norm); here $p_i=e_i\otimes e_i=<e_i, \cdot>e_i$ is the orthogonal projection onto the real line generated by $e_i$.
\item $(A^D)^D=A^D$ and $tr\left(A^DA\right)=tr((A^D)^2)$.
\item $tr(A^DB)=tr(AB)$ if $B$ is diagonal.
\end{enumerate}
\item The scalar manifold ${\Lambda}=\{\lambda 1:\;\lambda\in \mathbb R_{>0}\}$ is convex and closed in $\Sigma$, with tangent space at any $\lambda\in\Lambda$ given by $\mathbb R 1\subset {\cal H}_{\mathbb R}$.
\item A vector  $v=\mu+u$ is contained in $T_{d_0}{\Delta}^{\perp}$ if and only if $\mu=0$ and $u^D=0$. This follows from Remark \ref{tangente}, the fact that $\mu+u^D\in T_{d_0}\Delta$, and Remark (3) of this list. In other words for any $d_0\in\Delta$,
$$
T_{d_0}\Delta^{\perp}=T_1\Delta^{\perp}=\{v\in {\sf HS}^h:\; v \mbox{ is  codiagonal } \}=: \Gamma.
$$
\end{enumerate}

\begin{teo}\label{main2}Let $a\in{\sf HS}^h$. Then there exist $\lambda\in\mathbb R_{>0}$, $d\in \Delta$ and $x\in{\sf HS}^h$ such that:
$$
a+\lambda=(d+\lambda)\mbox{\rm e}^{(d+\lambda)^{-1}v}=(d+\lambda)^{\frac12} \mbox{\rm e}^{(d+\lambda)^{-\frac12}v(d+\lambda)^{-\frac12}}\;(d+\lambda)^{\frac12}.
$$
Moreover, for fixed $\lambda$, $d$ and $v$ are unique and $a+\lambda\mapsto (d,v)$ (which maps $\;\Sigma\to N{\Delta})$ is a real analytic isomorphism between manifolds.\end{teo}
\begin{proof}Let $\lambda=\|a\|_{\infty}+\epsilon$, for any $\epsilon>0$. Then $p=a+\lambda\in\Sigma$. Let $\Pi_{\Delta}(p)=d+\alpha$, where $d\in \Delta$. Now pick the unique $v\in T_{d+\alpha}{{\Delta}}^{\perp}$ such that $\mbox{\rm Exp}_{d+\alpha}(v)=p$, this operator $v$ has the desired form because of Remark (5) above. As a consequence of Lemma \ref{ef} $\alpha=\lambda$, for in this case $\beta=0$. \end{proof}

\smallskip

This theorem can be rephrased saying that, given a self-adjoint Hilbert-Schmidt operator $a$, for any $\lambda\in\mathbb R_{>0}$ such that $a+\lambda>0$, one has a unique factorization $a+\lambda =D\;{\rm e}^wD$ where $D=(\lambda+d)^{\frac12}>0$ is a diagonal operator and $w=D^{-1}vD^{-1}\in \Gamma$ is a self-adjoint operator with null diagonal. The normal bundle clearly splits in this case, so

\begin{prop} Consider the submanifolds $\Delta,\;\exp(\Gamma)\subset\Sigma$. Then the projection map $\Pi_{\Delta}$ induces a diffeomorphism $\Sigma\simeq \Delta\times \exp(\Gamma)$.
\end{prop}

\begin{coro}
For any $g\in {   {\cal H}_{\mathbb C}  }^{\bullet}$, there is a unique factorization $\;g=d{\rm e}^wu$, where $d$ is a positive invertible diagonal operator of ${\cal H}_{\mathbb C}$, $w$ is a self-adjoint operator with null diagonal in ${\cal H}_{\mathbb C}$ and $u$ is a unitary operator of ${\cal H}_{\mathbb C}$.
\end{coro}
\begin{proof}The previous results together with Theorem \ref{mpolar}.
\end{proof}

\section{A foliation of codimension one}\label{folia}

In this section we describe a foliation of the total manifold, and show how to translate the results from previous sections to a particular leaf (the submanifold $\Sigma_{1}$) in order to show an aplication concerning (finite dimensional) matrix algebras. Recall that we write ${{\sf HS}}^h$ for the self-adjoint Hilbert-Schmidt operators. Fix $\lambda\in\mathbb R_{>0}$. Let
$$
\Sigma_{\lambda}=\{a+\lambda\in \Sigma,\; a\in {{\sf HS}}^h \}.
$$
Observe that $\Sigma_{\lambda}\cap \Sigma_{\beta}=\emptyset$ when $\lambda\ne \beta$, since $a+\lambda=b+\beta$ implies $a-b=\beta-\lambda$. In this way, we can decompose the total space by means of these leaves, $\displaystyle\Sigma=\dot{\mathop{\cup}_{\lambda>0}}\Sigma_{\lambda}$.

\begin{prop} The leaves $\Sigma_{\lambda}$ are geodesically convex closed submanifolds.
\end{prop}
\begin{proof} We consider the projection to the convex scalar manifold $\Lambda$ (see Remark (4) above). The fact that the projection $\Pi_{{\Lambda}}$ is a contraction (therefore a continuous map) implies that $\Sigma_{\lambda}$ is closed; one must only observe that $\Sigma_{\lambda}=\Pi_{{\Lambda}}^{-1}(\lambda)$. To show that $\Sigma_{\lambda}$ is geodesically convex we recall that, by virtue of Lemma \ref{ef}, for any real $\lambda>0$ and any $p\in \Sigma_{\lambda}$, there is an identification \emph{via} the inverse exponential map at $p$, $T_p\Sigma_{\lambda}={{\sf HS}}^h$.\end{proof}

\begin{rem}\label{perpe}
Take $\delta+c\in {T_{a+\lambda}\Sigma_{\lambda}}^{\perp}$. Since ${T_{a+\lambda}\Sigma_{\lambda}}$ can be identified with ${{\sf HS}}^h$, the equality $\displaystyle
<\delta+c,d>_{_{a+\lambda}}=0\; \forall\; d \in {{\sf HS}}^h$ is equivalent to 
$$
tr\left[  (a+\lambda)^{-1}\left[(\delta+c)(a+\lambda)^{-1}-\delta/\lambda\right]d\right]=0\qquad \forall\; d \in {{\sf HS}}^h.
$$
Equivalently, $T_{a+\lambda}\Sigma_{\lambda}^{\perp}=span(a+\lambda)$; shortly $T_p\Sigma_{\lambda}^{\perp}=span(p)$ for any $p\in \Sigma_{\lambda}$.\end{rem}

\begin{prop}
Fix real $\alpha,\lambda>0$. Let $\Pi_{\alpha,\lambda}=\Pi_{\Sigma_{\lambda}}\left|_{\Sigma_{\alpha}}\right. :\Sigma_{\alpha}\to \Sigma_{\lambda}$. Then 
\begin{enumerate}
\item $\Pi_{\alpha,\lambda}(p)=\frac{\lambda}{\alpha}p$, so $\Pi_{\alpha,\lambda}(p)$ commutes with $p$.
\item  $\Pi_{\alpha,\lambda}$ is an isometric bijection between $\Sigma_{\alpha}$ and $\Sigma_{\lambda}$, with inverse $\Pi_{\lambda,\alpha}$. 
\item $\Pi_{\alpha,\lambda}$  gives parallel translation along vertical geodesics joining both leaves (that is, geodesics orthogonal to both leaves).
\end{enumerate}
\end{prop}
\begin{proof}
Notice that for a point $b+\alpha\in\Sigma_{\alpha}$ to be the endpoint of the geodesic $\gamma$, starting at $a+\lambda\in\Sigma_{\lambda}$, such that $L(\gamma)=\mbox{\rm dist}\left(b+\alpha,\Sigma_{\lambda}\right)$, we must have
$$
b+\alpha=\mbox{\rm Exp}_{a+\lambda}(x+c)=\mbox{\rm Exp}_{a+\lambda}(k.(a+\lambda))={\rm e}^k(a+\lambda)
$$
where $k\in\mathbb R$ comes from Remark \ref{perpe} above, since $x+c\in {T_{a+\lambda}\Sigma_{\lambda}}^{\perp}$. From Lemma \ref{ef}, we deduce that $k=\ln\left(\frac{\alpha}{\lambda}\right)$, and $a=\frac{\lambda}{\alpha} b$. So, $b+\alpha=\frac{\alpha}{\lambda}(a+\lambda)$ and also $\displaystyle \gamma(t)=(a+\lambda)\left(\textstyle\frac{\alpha}{\lambda}\right)^t$. Now it is clear that $\Pi_{\lambda}(b+\alpha)=\frac{\lambda}{\alpha}(b+\alpha)$ commutes with $b+\alpha$. To prove that $\Pi$ is isometric, observe that 
$$
\mbox{dist}(\Pi_{\alpha,\lambda}(p), \Pi_{\alpha,\lambda}(q))=\| \ln  ( (\frac{\lambda}{\alpha}p)^{-\frac12} (\frac{\lambda}{\alpha}q)  (\frac{\lambda}{\alpha}p)^{-\frac12} )     \|_{_2}=\| \ln  ( p^{-\frac12} q p^{-\frac12}  )     \|_{_2}=\mbox{dist}(p,q).
$$
That $\Pi$ gives parallel translation along $\gamma$ follows from the formula for $\Pi$ given in the first item of this proposition and the formula for the parallel translation given in Remark \ref{transporte}.\end{proof}

\smallskip

The normal bundle in the case of $M=\Sigma_1$ can be thought of as a direct product:

\begin{prop}\label{isometrica}
The map $\,T:\Sigma\to \Sigma_1\times \Lambda\,$, which assigns $a+\alpha\mapsto\left(\frac{1}{\alpha}(a+\alpha) ,\alpha \right)\;$ is bijective and isometric ($\Sigma_1$ and $\Lambda$ have the induced submanifold metric). In other words, there is a Riemannian isomorphism $\Sigma\simeq \Sigma_1\times \Lambda$.
\end{prop}

\begin{prop}\label{paralelas}
The leaves $\Sigma_{\alpha}$, $\Sigma_{\lambda}$ are also parallel in the following sense: any minimizing geodesic joining a point in one of them with its projection in the other is orthogonal to both of them. See Figure 1 below.

\begin{figure}[ht]
\centerline{\includegraphics[height=3.845cm,width=6.15cm]{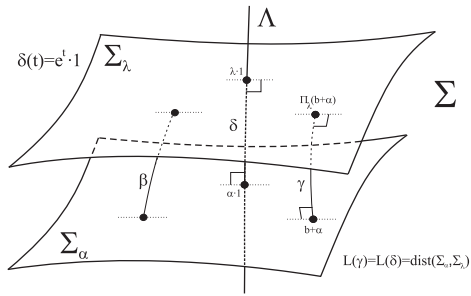}}
\caption{\small The geodesics $\gamma$ and $\delta$ are minimizing, the geodesic $\beta$ is not.}
\end{figure}

\smallskip

For any $b+\alpha\in \Sigma_{\alpha}$ we have $\displaystyle \mbox{\rm dist}(b+\alpha,\Sigma_{\lambda})= \mbox{\rm dist}(\Sigma_{\alpha},\Sigma_{\lambda})=\mid\ln\left(\textstyle\frac{\alpha}{\lambda}\right)\mid$. In particular, the distance between $\alpha,\lambda$ in the scalar manifold $\Lambda$ is given by the Haar measure  of the open interval $(\alpha,\beta)$  on the multiplicative group $\mathbb R_{>0}$.
\end{prop}
\begin{proof} It is a straightforward computation that follows from the previous results; the last statement was observed by E. Vesentini in another context \cite{ves}.\end{proof}

\smallskip

Since $\Sigma$ is a symmetric space, curvature is preserved when we parallel-translate bidimensional planes. Note also that vertical planes (\emph{i.e.} planes generated by a vector $v\in {\sf HS}^h=T_{\lambda}\Sigma_{\lambda}$ and $\lambda$) are commuting sets of operators.

\begin{prop}\label{flatso}
Let $p\in\Sigma_{\lambda}$. Then the sectional curvature of vertical 2-planes is zero.
\end{prop}
\begin{proof}
It follows from the formula for the curvature given in Section \ref{ctensor}.
\end{proof}

\subsection{The embedding of $M_n^+$ in $\Sigma_1$}\label{finite}

Let $M_n^+$ be the set of positive invertible $n\times n$ matrices (see the introduction of this paper). First note that we can embed $M_n^+\hookrightarrow\Sigma_1$ for any $n\in\mathbb N$: fix an orthonormal basis $\{e_n\}_{n\in\mathbb N}$ of $H$, let $p_{ij}=e_i\otimes e_j$, and identify the set $M_n$ of real $n\times n$ matrices with the set 
$$
{\cal T}=\{\sum_{i,j=1}^n a_{_{ij}}\; p_{ij} \;: \;a_{_{ij}}=a_{_{j i}}\in\mathbb R\}\subset {{\sf HS}}^h.
$$
We identify the manifold $M_n^+$ with ${\cal P}=\left\{\mbox{\rm e}^T : T\in {\cal T}\right\}\subset \Sigma_1
$ and the tangent space at each $\mbox{\rm e}^T\in{\cal P}$ is ${\cal T}$. The set ${\cal P}$ is closed and convex in $\Sigma_1$ by Corollary \ref{triple}. Let us call $S={\rm span}(e_1,\cdots,e_n)$, $S^{\perp}={\rm span}(e_{n+1},e_{n+2}\cdots )$. The operator $P_{_S}$ is the orthogonal projection to $S$ and $Q_{_S}=1-P_{_S}$ is the orthogonal projection to $S^{\perp}$. Using matrix blocks, for any operator $A\in {\rm B}(S)$, we identify
$$ 
\,{\cal T}=\left\{\left( \begin{array}{lc} A &  0 \\  0 & 0 \end{array} \right)\right\}\quad\mbox{ and }\quad{\cal P}=\left\{\left( \begin{array}{lc} \mbox{\rm e}^A &  0 \\  0 & 1 \end{array} \right)\right\}.
$$

\begin{rem}\label{ortog} There is a direct sum decomposition of ${{\sf HS}}^h=\,{\cal T}\oplus {\cal J}$
where operators in $J \in {\cal J}$ are such  that $P_{_S}\,J\,P_{_S}=0$.  A straightforward computation using the matrix-block representation shows that $<a,b>_{_2}=0$ for any $a\in {\cal T}, b\in {\cal J}$, which says ${\cal T}^{\perp}={\cal J}$ (here we consider ${{\sf HS}}^h$ as the total space). So the manifolds ${\rm exp}({\cal J})$ and ${\cal P}={\rm exp}({\cal T})$ are orthogonal at $1$, the unique intersection point. In the notation of Theorem \ref{homo}, it is also clear that ${\cal P}=M_G\simeq G/K$, where
$$
G=\left( \begin{array}{cc}  GL(n, \mathbb C) &  0 \\  0 & 1 \end{array} \right)\quad\mbox{and }\; K=\left( \begin{array}{cc} {\cal U}(n,\mathbb C) &  0 \\  0 & 1 \end{array} \right).
$$
\end{rem}

\begin{teo}\label{mn}
Let ${\cal P}\simeq M_n^+\subset\Sigma_1$ with the above identification. Then for any positive invertible operator $\mbox{\rm e}^b\in\Sigma_1$, ($b\in {{\sf HS}}^h$) there is a unique factorization of the form
$$ 
\mbox{\rm e}^{b}=\left( \begin{array}{lc} \mbox{\rm e}^A &  0 \\  0 & 1 \end{array} \right)
\;\mbox{\rm exp}  \left\{    \left( \begin{array}{lc} \mbox{\rm e}^{-A} &  0 \\  0 & 1 \end{array} \right)   \left( \begin{array}{lc} 0 &  Y^* \\  Y & X \end{array} \right)  \right\}\quad \mbox{ where if } a=\left( \begin{array}{lc} A &  0 \\  0 & 0 \end{array} \right)\in {\cal T}
$$
then $\mbox{\rm e}^a=\mbox{\rm e}^AP_{_S}+Q_{_S}\in {\cal P}\simeq M_n^+$, $X^*=X$ is a Hilbert-Schmidt opertor acting on the Hilbert space $S^{\perp}$, and $Y\in {\rm B}(S,S^{\perp})$. 
\end{teo}

An equivalent expression for the factorization is
$$ 
\mbox{\rm e}^{ b}=\left( \begin{array}{lc} \mbox{\rm e}^{A/2} &  0 \\  0 & 1 \end{array} \right)
\;\mbox{\rm exp}\left\{       \left( \begin{array}{lc}  0&  \mbox{\rm e}^{-A/2} Y^* \\  Y \mbox{\rm e}^{-A/2}  & X \end{array} \right)  \right\} \left( \begin{array}{lc} \mbox{\rm e}^{A/2} &  0 \\  0 & 1 \end{array} \right).
$$
\begin{proof}
From previous theorems and the observations we made, we know that
${\rm e}^{ b}= {\rm e}^{a/2} C {\rm e}^{a/2}$, where
$$
C=\mbox{\rm exp}\left\{    \left( \begin{array}{lc} \mbox{\rm e}^{-A/2} &  0 \\  0 & 1 \end{array} \right)
   \left( \begin{array}{lc} V_{11} &  V_{21}^* \\  V_{21}  & V_{22} \end{array} \right) \left( \begin{array}{lc} \mbox{\rm e}^{-A/2} &  0 \\  0 & 1 \end{array} \right) \right\} 
$$
for some $A\in {\rm B}(S)$ and some $v\in {{\sf HS}}^h$. That $V_{11}=0$ follows from the fact (see Remark \ref{ortog}) that ${\cal T}^{\perp}={\cal J}$, and $v\in T_{{\rm e}^a}{\cal P}^{\perp}$ if and only if $tr({\rm e}^{-A}B{\rm e}^{-A}V_{11})=0$ for any $B\in {\cal T}$. \end{proof}

\begin{rem}For any $b\in {{\sf HS}}^h$, the operator
$$ 
\mbox{\rm e}^{ a}=\mbox{\rm e}^AP_{_{S}^{ \;}}+P_{S^{\;^\perp}}=\left( \begin{array}{lc} \mbox{\rm e}^{A} &  0 \\  0 & 1 \end{array} \right)=\mbox{\rm exp}\left( \begin{array}{lc} A &  0 \\  0 & 0 \end{array} \right)
$$
is the 'first block' $n\times n$ matrix which is closest to ${\rm e}^b$ in $\Sigma$, and with a slight abuse of notation for the traces of ${\rm B}(S)$ and ${\rm B}(S^{\perp})$, we have
$$
{\mbox{\rm dist}({\cal P},{\rm e}^b)=\mbox{\rm dist}(\mbox{\rm e}^a,\mbox{\rm e}^b) = \sqrt{\|Y\,{\rm e}^{-A/2}\|_{_2}^{^2}+\|X \|_{_2}^{^2}}}.
$$
\end{rem}

\begin{coro}
For any $g\in {   {\cal H}_{\mathbb C}  }^{\bullet}$ there is a unique factorization $g=\lambda r {\rm e}^v u$, where $\lambda \in\mathbb R_{>0}$, $u\in\cal U({\cal H}_{\mathbb C})$ is a unitary operator,
$$
r=\left( \begin{array}{lc} R &  0 \\  0 & 1 \end{array} \right)\qquad\qquad v=\left( \begin{array}{lc} 0 & Y^* \\  Y & X \end{array} \right)
$$
with $R\in {\rm B}(S)^+\simeq M_n^+$, $X=X^*\in {\rm B}(S^{\perp})$ a Hilbert Schmidt, and $Y\in {\rm B}(S,S^{\perp}).$
\end{coro}
\begin{proof}
We use the notation of Remark \ref{ortog}. Note that, by Theorem \ref{mpolar}, $g=r {\rm e}^x u$ with $u\in\cal U({\cal H}_{\mathbb C})$, $r\in {\cal P}=\exp({\cal T})$ and $x\in{\cal T}^{\perp}$. But ${\cal T}^{\perp}={\cal J}\oplus \mathbb R1$ if we consider ${\cal H}_{\mathbb R}$ as the total space, and ${\rm e}^{\alpha+a}={\rm e}^{\alpha}{\rm e}^a$ if $\alpha\in\mathbb R$.
\end{proof}

\vspace*{.5cm}

{\sc\small\noindent Gabriel Larotonda\newline
Instituto de Ciencias, Universidad Nacional de General Sarmiento.\newline
JM Gutiérrez 1150 (1613) Los Polvorines. Buenos Aires, Argentina.}

\smallskip

{\sf\small \noindent e-mail: glaroton@ungs.edu.ar\newline
Tel/Fax: (+54-011)-44697501}


\medskip


\begin{thebibliography}{XX}

\label{bib}

\bibitem{acms} E. Andruchow, G. Corach, M. Milman and D. Stojanoff, {\it Geodesics and interpolation}, Revista de la Unión Matemática Argentina \textbf{40} (1997) n$^o$3 and 4, 83-91. 

\bibitem{gandru} E. Andruchow and G. Larotonda, {\it Nonpositively Curved Metric in the Positive Cone of a Finite von Neumann Algebra}, J. London Math. Soc., to appear.

\bibitem{atkin1} C.J. Atkin, {\it The Hopf-Rinow theorem is false in infinite dimensions}, Bull. London Math. Soc.  (1975) n$^o$7, 261-266.

\bibitem{atkin2} C.J. Atkin, {\it Geodesic and metric completeness in infinite dimensions}, Hokkaido Math. J. \textbf{26} (1997), 1-61.

\bibitem{ball0} W. Ballmann, {\it Spaces of Nonpositive Curvature}, Jahresber. Deutsch. Math. Verein. \textbf{103} (2001), n$^o$2, 52-65.

\bibitem{bhatia} R. Bhatia, {\it On the exponential metric increasing property}, Linear Algebra Appl. 375 (2003), 211-220.

\bibitem{caldera} P. Calderón, {\it Intermediate spaces and interpolation, the complex method}, Studia Math. \textbf{24} (1964), 113-190.

\bibitem{cm1} G. Corach and A. Maestripieri, {\it Differential and metrical structure of positive operators}, Positivity \textbf{3} (1999) 297-315.

\bibitem{cm2} G. Corach and A. Maestripieri, {\it Positive operators on Hilbert space: a geometrical view point}, Colloquium on Homology and Representation Theory, Bol. Acad. Nac. Cienc. (Córdoba) \textbf{65} (2000) 81-94.  

\bibitem{cpr6} G. Corach, H. Porta and L. Recht, {\it Differential Geometry of Spaces of Relatively Regular Operators}, Integral Equations Operator Theory (1990) n$^o$13, 771-794.

\bibitem{cpr7} G. Corach, H. Porta and L. Recht, {\it Splitting of the Positive Set of a C$^*$-Algebra}, Indag. Math. NS 2 (1991) n$^o$4, 461-468.  

\bibitem{segal} G. Corach, H. Porta and L. Recht, {\it A Geometric interpretation of Segal's inequality $\|{\rm e}^{x+y}\|\le\|{\rm e}^{x/2}{\rm e}^y{\rm e}^{x/2}\|$}, Proc. of the AMS \textbf{115} (1992) n$^o$1, 229-231

\bibitem{cpr2} G. Corach , H. Porta and L. Recht, {\it The Geometry of the Space of Selfadjoint Invertible Elements in a C$^*$-algebra}, Integral Equations Operator Theory \textbf{16} (1993), 333-359.

\bibitem{ebe2} P. Eberlein, {\it Geometry of Nonpositively Curved Manifolds}, Chicago Lectures in Mathematics, University of Chicago Press, Chicago, IL  (1996).

\bibitem{grossman1} N. Grossman, {\it Hilbert manifolds without epiconjugate points}, Proc. AMS \textbf{16} (1965), 1365-1371.

\bibitem{sigurdur} S. Helgason, {\it Differential Geometry, Lie Groups and Symmetric Spaces}, Academic Press, New York (1962).

\bibitem{pharpe} P. de la Harpe, {\it Classical Banach-Lie Algebras and Banach-Lie Groups of Operators in Hilbert Space}, Lecture Notes in Mathematics \textbf{285}, Springer, Berlin (1972). 

\bibitem{jost} H. Jost, {\it Nonpositive Curvature: Geometric and Analytic Aspects}, Lectures in Mathematics, Birkhäuser, Berlin (1997).

\bibitem{lang0} S. Lang, {\it Introduction to differentiable manifolds}, Interscience, New York (1962).

\bibitem{lang} S. Lang, {\it Differential and Riemannian Manifolds}, Springer-Verlag, Berlin-New York (1995).

\bibitem{lawson} J. Lawson and Y. Lim, {\it Symmetric spaces with convex metrics}, preprint (2006).

\bibitem{mcalpin} J. McAlpin, {\it Infinite Dimensional Manifolds and Morse Theory}, PhD. Thesis, Columbia University (1965).

\bibitem{mostow} G.D. Mostow, {\it Some new decomposition theorems for semi-simple groups}, Mem. Amer. Math. Soc. \textbf{14} (1955), 31-54.

\bibitem{neeb} K.H. Neeb, {\it A Cartan-Hadamard Theorem for Banach-Finsler manifolds}, Geom. Ded. \textbf{95} (2002), 115-150.

\bibitem{pr1} H. Porta and L. Recht, {\it Spaces of Projections in Banach Algebras}, Acta Matemática Venezolana \textbf{38} (1987), 408-426.

\bibitem{pr4} H. Porta and L. Recht, {\it Conditional Expectations and Operator Decompositions}, Ann. Global Anal. Geom. \textbf{12} (1994), 335-339.

\bibitem{simon} B. Simon, {\it Trace ideals and their applications}, London Mathematical Society Lecture Note Series, 35. Cambridge University Press, Cambridge-New York (1979).

\bibitem{ves} E. Vesentini, {\it Invariant metrics on convex cones}, Ann. Scuola Norm. Sup. Pisa Cl. Sci. (Ser.4) \textbf{3} (1976), 671-696.

\end{thebibliography}
\end{document}